\documentclass[10pt,journal,final]{IEEEtran}
\usepackage{amssymb,amsmath,amsthm}
\usepackage{graphicx}
\usepackage{enumerate}
\usepackage{comment,cite,color}
\usepackage{cite,color}
\usepackage{mathrsfs}
\usepackage{epsfig}
\usepackage{lscape}
\usepackage{subfigure}
\usepackage{epstopdf}
\usepackage{caption}
\usepackage{array}
\usepackage{algorithm}
\usepackage{algpseudocode}
\usepackage{bm}
\usepackage{booktabs}
\usepackage{threeparttable}
\usepackage{multirow}
\usepackage{bbm}
\begin{document}
% =========================Def===============================
	
	\newcommand{\PP}{\mathbb{P}}
	\newcommand{\HH}{\mathbf{H}}
	\newcommand{\innerp}[1]{\langle {#1} \rangle}
	\newcommand{\norm}[1]{\|{#1}\|_2}
	\newcommand{\abs}[1]{\lvert#1\rvert}
	\newcommand{\absinn}[1]{\vert\langle {#1} \rangle\rvert}
	\newcommand{\argmin}[1]{\mathop{\rm argmin}\limits_{#1}}
	\newcommand{\argmax}[1]{\mathop{\rm argmax}\limits_{#1}}
	\newcommand{\ds}{\displaystyle}
	\renewcommand{\Re}{{\rm Re}}
	\renewcommand{\Im}{{\rm Im}}
	\newcommand{\wt}{\widetilde}
	\newcommand{\ra}{{\rightarrow}}
	\newcommand{\lra}{{\longrightarrow}}
	\newcommand{\eproof}{\hfill\rule{2.2mm}{3.0mm}}
	\newcommand{\esubproof}{\hfill$\Box$}
	\newcommand{\Proof}{\noindent {\bf Proof.~~}}
	\newcommand{\D}{{\mathcal D}}
	\newcommand{\TT}{{\mathcal T}}
	\renewcommand{\SS}{{\mathcal S}}
	\newcommand{\BS}{{\mathbb S}}
	\newcommand{\Prob}{{\rm Prob}}
	\newcommand{\R}{{\mathbb R}}
	\newcommand{\0}{{\mathbf 0}}
	\newcommand{\Z}{{\mathbb Z}}
	\newcommand{\T}{{\mathbb T}}
	\newcommand{\C}{{\mathbb C}}
	\newcommand{\CG}{{\mathcal G}}
	\newcommand{\Q}{{\mathbb Q}}
	\newcommand{\N}{{\mathbb N}}
	\newcommand{\ep}{\varepsilon}
	\newcommand{\wmod}[1]{\mbox{~(mod~$#1$)}}
	\renewcommand{\eqref}[1]{(\ref{#1})}
	\newcommand{\inner}[1]{\langle #1 \rangle}
	\newcommand{\shsp}{\hspace{1em}}
	\newcommand{\mhsp}{\hspace{2em}}
	\newcommand{\FT}[1]{\widehat{#1}}
	\newcommand{\conj}[1]{\overline{#1}}
	\newcommand{\ber}{\nu_{\lambda}}
	\newcommand{\bern}{\nu_{\lambda_n}}
	\newcommand{\biasber}{\nu_{\lambda, p}}
	\newcommand{\bequiv}{\sim_\lambda}
	\newcommand{\bnequiv}{\sim_{\lambda_n}}
	
	\newcommand{\E}{{\mathbb E}}
	\newcommand{\B}{{\mathcal B}}
	\newcommand{\MC}{{\mathcal C}}
	\newcommand{\ML}{{\mathcal L}}
	\newcommand{\MLP}{{\mathcal{PL}}}
	\newcommand{\MLS}{{\mathcal{PS}}}
	\newcommand{\NN}{{\mathcal N}}
	\renewcommand{\i}{{\mathbf i}}
	\newcommand{\rank}{{\rm rank}}
	\newcommand{\supp}{{\rm supp}}
	\newcommand{\diag}{{\rm diag}}
	\newcommand{\ph}{{\rm ph}}
	\renewcommand{\j}{{\mathbf j}}
	\newcommand{\vx}{{\mathbf x}}
	\newcommand{\vy}{{\mathbf y}}
	\newcommand{\vz}{{\mathbf z}}
	\newcommand{\vO}{{\mathbf 0}}
	\newcommand{\vo}{{\mathbf o}}
	\newcommand{\va}{{\mathbf a}}
	\newcommand{\vb}{{\mathbf b}}
	\newcommand{\vd}{{\mathbf d}}
	\newcommand{\vv}{{\mathbf v}}
	\newcommand{\vu}{{\mathbf u}}
	\newcommand{\vf}{{\mathbf f}}
	\newcommand{\vg}{{\mathbf g}}
	\newcommand{\vw}{{\mathbf w}}
	\newcommand{\ve}{{\mathbf e}}
	\newcommand{\vh}{{\mathbf h}}
	\newcommand{\vr}{{\mathbf r}}
	\newcommand{\vepsilon}{{\bm \epsilon}}
	\renewcommand{\H}{{\mathbb H}}
	\newcommand{\qHH}{\widetilde{\HH}}
	\newcommand{\qC}{\widetilde{\C}}
	\newcommand{\ul}[1]{\underline{#1}}
	\newcommand{\MM}{\mathbf M}
	\newcommand{\XX}{\mathbf X}
	\newcommand{\x}{\mathbf x}
	\newcommand{\ba}{\mathbf a}
	\newcommand{\F}{{\mathbb F}}
	\newcommand{\G}{{\mathcal G}}
	\newcommand{\CR}{{\mathcal R}}
	\newcommand{\sn}{{\mathbb S}^n}
	\newcommand{\vn}{{\mathbf v}^n}
	\newcommand{\un}{{\mathbf u}^n}
	\newcommand{\uD}{\underline{D}}
	\newcommand{\dd}{{\mathrm d}}
	\newcommand{\s}{{\mathbf s}}
	
	\renewcommand{\thefigure}{\arabic{figure}}
	\newtheorem{definition}{Definition}[section]
	\newtheorem{corollary}{Corollary}[section]
	\newtheorem{theorem}{Theorem}[section]
	\newtheorem{example}{Example}[section]
	\newtheorem{lemma}{Lemma}[section]
	\newtheorem{remark}{Remark}[section]
	\newtheorem{notation}{Notation}[section]
	\newtheorem{prop}{Proposition}[section]
	\newcommand{\zz}{^{\top}}
	\newcommand{\sign}{{\mathrm sign}}
	\newcommand{\red}[1]{\textcolor{red}{#1}}
	
	% =========================Def===============================
	
	% =========================Title=============================
		
	\title{Perturbed amplitude flow for phase retrieval
	%\thanks{Research of Zhiqiang Xu was supported  by NSFC grant ( 11422113,  91630203, 11331012) and by National Basic Research Program of China (973 Program 2015CB856000).}
	\thanks{B.~Gao is with School of Mathematical Sciences, Nankai University, Tianjin, China. Email: gaobing@nankai.edu.cn}
	\thanks{X.~Sun is with Microsoft Research Asia Email: xinsun@microsoft.com}
	\thanks{ Y.~Wang is with Department of Mathematics, The Hong Kong University of Science and Technology, Clear Water Bay, Kowloon, Hong Kong. Email:  yangwang@ust.hk}
	\thanks{ Z.~Xu is with Inst. Comp. Math., Academy of Mathematics and Systems Science, Chinese Academy of Sciences, Beijing, 100091, China; School of Mathematical Sciences, University of Chinese Academy of Sciences, Beijing 100049, China. Email:  xuzq@lsec.cc.ac.cn}			
	\thanks{Yang Wang was supported in part by the Hong Kong Research Grant Council grants 16306415 and 16308518. Zhiqiang Xu was supported  by NSFC grant (91630203, 11688101), Beijing Natural Science Foundation (Z180002).}
	}	
	\author{ Bing Gao,\,\,Xinwei Sun,\,\,Yang Wang,\,\, Zhiqiang Xu}
	\date{\today}
	\markboth{Date of current version Sep.~2020}%
	{Shell \MakeLowercase{\textit{et al.}}: Bare Demo of IEEEtran.cls for IEEE Journals}
    \maketitle
	% =========================Title=============================
	
	% =========================Abs===============================
	\begin{abstract}
	In this paper, we propose a new non-convex algorithm for solving the phase retrieval problem, i.e., the reconstruction of a signal $ \vx\in\H^n $ ($\H=\R$ or $\C$) from phaseless samples $ b_j=\abs{\langle \va_j, \vx\rangle } $, $ j=1,\ldots,m $. The proposed algorithm solves a new proposed model, perturbed amplitude-based model, for phase retrieval and is correspondingly named as {\em Perturbed Amplitude Flow} (PAF). We prove that PAF can recover $c\vx$ ($\abs{c} = 1$) under $\mathcal{O}(n)$ Gaussian random measurements (optimal order of measurements).
	Starting with a designed initial point, our PAF algorithm iteratively converges to the true solution at a linear rate for both real and complex signals. Besides, PAF algorithm needn't any truncation or re-weighted procedure, so it enjoys simplicity for implementation. The effectiveness and benefit of the proposed method are validated by both the simulation studies and the experiment of recovering natural images.	
	\end{abstract}
	% =========================Abs===============================
	
	% =========================Sec===============================
	\section{Introduction}
	
	\subsection{Problem Setup and Related Work}
	\IEEEPARstart{I}{n} this paper, we consider the well-known \emph{phase retrieval} problem, which aims to recover a signal $ \vx\in \H^n $, where $\H=\R$ or $\C$, from phaseless measurements 	
	\[
	b_j = \abs{\langle \va_j, \vx\rangle},\quad  j=1,\ldots, m.
	\]
	Here $\vx\in\H^n$ is the {\em target signal} or {\em the target vector} and the vectors $ \va_j \in\H^n$ for all $j$ are the {\em measurement vectors}. Phase retrieval has many applications in both science and engineering, such as X-ray crystallography \cite{miao1999extending, elser2018benchmark},  astronomy \cite{fienup1987phase}, optics \cite{walther1963question, millane1990phase}, microscopy \cite{miao2008extending}.
	
	Due to the removal of phase information in the measurements  $| \langle \va_j, \vx\rangle| $, we can only recover $\vx$ up to a unimodular constant. Moreover, it is also known that $ \mathcal{O}(n) $ general measurements are enough to recover a signal $\vx \in\H^n$ uniquely. Particularly, it was shown that  $m\geq 2n-1$ and $m\geq 4n-4$ generic measurements $\{\va_j\}_{j=1}^m\subset \H^n$ are sufficient to recover any $\vx\in\H^n$ up to a unimodular constant for $\H=\R$ and $\H=\C$, respectively  \cite{balan2006signal,bandeira2014saving, WX17}.
	
	The original phase retrieval problem mainly considers the recovery of a signal from its Fourier transform magnitude \cite{huang2016phase}  or the magnitude of the short-time Fourier transform
	\cite{eldar2015sparse, jaganathan2016stft, bendory2018non-convex}. At the same time, more algorithms have been developed for general cases, in which random observations are considered, which also provide heuristic algorithms for practical applications. They can be roughly divided into two categories: the convex methods and the non-convex ones. For convex methods, the general strategy is to lift the phase retrieval problem into a problem of recovering a rank-one matrix and apply the semi-definite programming to solve it. The first such method, called  PhaseLift \cite{candes2013phaselift,candes_phaselift,candes_matrix_completion}, can achieve the exact recovery using $ m=\mathcal{O}(n) $ independent  Gaussian random measurements $ \va_j $, $ j=1,\ldots,m $. However, such an approach is computationally inefficient for large dimensional problems since semi-definite programming for $ n\times n $ matrices is slow for large $n$.  An alternative method called PhaseMax \cite{goldstein2018phasemax,dhifallah2017fundamental,hand2016elementary}  aims to recover the signal $ \vx $ by solving the model
	\begin{equation}\label{phasemax_model}
	\max_\vz\,\, \Re(\langle \vz, \hat{\vz})) \quad \text{subject to} \quad \abs{\langle \va_j, \vz\rangle} \leq b_j,
	\end{equation}
	where $ \hat{\vz} $ is an approximation to the true signal $ \vx $. It is proved that this method can recover $\vx$ with high probability when $ m\geq 4n/\theta $ where $ \theta = 1-\frac{2}{\pi}\text{angle}(\langle\hat{\vz},\vx\rangle) $. However, numerical experiments have shown that larger oversampling ratios $ m/n $ are often required for exact recovery, especially compared to several non-convex algorithms.
	
	In a different direction, a series of non-convex approaches have been proposed and studied. Among such schemes, early studies are based on the alternating projection approach, including the works by Gerchberg and Saxton \cite{gerchberg1972practical} and Fineup \cite{fienup1982phase}. These methods often perform well numerically but lack theoretical foundations. Motivated by the success of alternating minimization, Netrapalli et al \cite{netrapalli2013phase} developed the AltMinPhase method that is shown to achieve linear convergence with $\mathcal{O}(n\log^3 n) $ Gaussian random  measurements and resampling. Recently, the sample complexity is improved to $\mathcal{O}(n)$ Gaussian random measurements in complex number field under a carefully chosen initial point by Waldspurger in \cite{waldspurger2018phase}. However, such an alternating projection-based approach also suffers from larger computational complexity, due to the projection step.	
	More recently another framework was proposed, in which one starts from a ``good'' initial guess and try to iteratively refine it by solving a given model such as the intensity-based model \cite{candes_wf, gao2017phaseless}
	\begin{equation}\label{wf_model}
	\min_{\vz}\,\,g(\vz):= \frac{1}{4m}\sum_{j=1}^{m}\left(|\va_j^*\vz|^2 -b_j^2 \right)^2,
	\end{equation}
	or the amplitude-based model \cite{wang2018solving,zhang2016reshaped,huangxu,wang2018phase}
	\begin{equation}\label{amplitude_model}
	\min_{\vz}\,\,f(\vz):= \frac{1}{2m}\sum_{j=1}^{m}\left(|\va_j^*\vz| -b_j \right)^2,
	\end{equation}
	or the Poisson likelihood model \cite{chen2015solving}
	\begin{equation}\label{likelihood_model}
	\min_{\vz}\,\,h(\vz):=-\sum_{j=1}^{m} \left( b_j^2\log(|\va_j^*\vz|^2)-|\va_j^*\vz|^2\right).
	\end{equation}
	
	Among existing proposed algorithms to solve \textbf{intensity-based model} (\ref{wf_model}), Cand\`{e}s et al have developed the Wirtinger Flow method (WF) \cite{candes_wf} to recover $ \vx $ via gradient descent. It achieved provable linear convergence with $ m=\mathcal{O}(n\log n)$ Gaussian random measurements under carefully chosen initialization method.
	Particularly, Sun et al \cite{sun2018geometric} proved the benign geometric landscape of~\eqref{wf_model} under $\mathcal{O}(n \mathrm{poly}(\log{n}))$ Gaussian measurements, motivating the Trust-region method to avoid spurious local minimizers. Besides, Ma et al \cite{2017Implicit} proved the ``nice'' geometry of~\eqref{wf_model} under Gaussian random measurements, explaining the favorable performance of unregularized gradient descent. Such geometric benefits guarantee the success of gradient descent for this non-convex phase retrieval problem. 
	Recently, the result is refined in real field $\H=\R$ to achieve a reduction of measurements $m = \mathcal{O}(n) $ by solving \textbf{amplitude-based model} (\ref{amplitude_model}) via gradient descent \cite{zhang2016reshaped} or via truncated gradient descent \cite{wang2018solving} or via reweighted gradient descent
	\cite{wang2018phase}, or by solving \textbf{Poisson likelihood model} (\ref{likelihood_model}) via modified gradient descent \cite{chen2015solving}.  In detail, Zhang et al \cite{Zhang2017A} have proposed Reshaped Wirtinger Flow,  which named  Amplitude Flow (AF) in this paper to coincide with the model used,  to solve model (\ref{amplitude_model}) by gradient descent.  Wang et al \cite{wang2018solving} have proposed Truncated Amplitude Flow (TAF) to solve model (\ref{amplitude_model}) by truncated gradient descent. Wang et al \cite{wang2018phase} have designed Reweighted Amplitude Flow (RAF) to solve model (\ref{amplitude_model}) via reweighted gradient descent. Chen and Cand\`{e}s \cite{chen2015solving} have designed Truncated Wirtinger Flow (TWF), which solves model (\ref{likelihood_model}) by modified gradient descent.

	From the perspective of theoretical analysis, the methods that given in AF, TAF, RAF and TWF all can achieve linear convergence under the optimal order of measurements. Different from truncation-based methods (e.g., TAF \cite{wang2018solving}, TWF \cite{chen2015solving}) that remove the components having too much influence on the search direction, the RAF \cite{wang2018phase} implements re-weighted procedure to control such components by reducing their weights at each update. Instead of using truncation or re-weighted procedures to get reliable gradients, the AF \cite{zhang2016reshaped} method performs gradient descent directly.  But the analysis of AF is based on the fact that the value of $ \text{sign}(\langle\va_j,\vz\rangle) $ equals $ -1 $ or $ 1 $, which can only be satisfied in the real number field. This fact is also required by TAF, TWF and RAF. Thus the theoretical results can't be extended to the complex case trivially.
	
	In this paper, we introduce a new perturbed amplitude-based
	model to address these theoretical deficiencies and limitations in this framework.	
%%-------------------

	\subsection{Our Contribution: The Perturbed Amplitude Flow (PAF)}
	
	We propose the {\em Perturbed Amplitude Flow (PAF)} algorithm in this paper through the following model:
	\begin{equation}\label{model}
	\min_{\vz}f_{\vepsilon}(\vz):= \min_{\vz} \frac{1}{m}\sum_{j=1}^{m}\left(\sqrt{|\va_j^*\vz|^2+\epsilon_j^2} - \sqrt{ b_j^2+\epsilon_j^2} \right)^2,
	\end{equation}
	where $ \vepsilon=[\epsilon_1,\ldots,\epsilon_m]\in\R^m $ have prescribed value, with the requirement that
	\begin{equation}\label{condition_eps}
	\epsilon_j\neq 0\quad \text{for all} \quad b_j\neq 0.
	\end{equation}
	Note that if $b_j=0$, then $\left(\sqrt{|\va_j^*\vz|^2+\epsilon_j^2} - \sqrt{ b_j^2+\epsilon_j^2} \right)^2$ is smooth regardless of the value of $\epsilon_j$, even when $\epsilon_j=0$. The loss function $ f_{\vepsilon}$ is thus smooth. When all $\epsilon_j=0$, this model is reduced to the classic amplitude-based model (\ref{amplitude_model}). So we shall name it as the {\it perturbed amplitude-based model} and name the corresponding gradient descent method as Perturbed Amplitude Flow (PAF).
	
	In the perturbed amplitude-based model (\ref{model}), $ \vepsilon $ not only keeps the loss function smooth but plays a role similar to {\it truncation/re-weighted} while reducing the effects of bad observations. From the previous work \cite{wang2018solving, wang2018phase},
	we know that only the gradients associated with sizable $ |\va_j^*\vz|/|\va_j^*\vx| $ offer meaningful directions. In detail, considering the model (\ref{amplitude_model}), when $ \H=\R $ the wirtinger derivative of $ f $  concerning to $ \vz $ is
	\begin{align*}
	\nabla f(\vz) &= \frac{1}{m}\sum_{j=1}^{m}\left(1-\frac{b_j}{|\va_j^*\vz| }\right)\va_j\va_j^*\vz\\
	&=\frac{1}{m}\sum_{j=1}^{m}\left(\va_j\va_j\zz\vh+\va_j|\va_j\zz\vx|\Big(\frac{\va_j\zz\vx}{|\va_j\zz\vx|}-\frac{\va_j\zz\vz}{|\va_j\zz\vz|}\Big)\right),
	\end{align*}
	with $ \vh = \vz-\vx $. Note that the first term $ \va_j\va_j\zz\vh
	 $ flows a desirable direction, whereas the second term $ \va_j|\va_j\zz\vx|\Big(\frac{\va_j\zz\vx}{|\va_j\zz\vx|}-\frac{\va_j\zz\vz}{|\va_j\zz\vz|}\Big) $ has negative influence and such an influence can be reduced when $ \va_j\zz\vz $ shares the same sign with $ \va_j\zz\vx $. The TAF \cite{wang2018solving} established that those terms with inconsistent sign are normally those terms with small $ |\va_j\zz\vz| $ in real case, which motivates a truncation scheme that drops the terms with small $ |\va_j\zz\vz|/|\va_j\zz\vx| $. Instead of abandoning those gradients, RAF \cite{wang2018phase} uses re-weighted procedure to reduce the influence of those components. However, these analyses heavily rely on the sign of each element equal 1 or -1, therefore hard to be extended to the complex case.
	
	For our model, with a suitable choice of $ \vepsilon $, one can control the size of the gradient. This is essential for avoiding the extremely large gradient components. More precisely, note that the Wirtinger derivative of $ f_{\vepsilon} $ with respect to $ \vz $ is
	\begin{align*}
	\nabla f_{\vepsilon}(\vz) = \frac{1}{m}\sum_{j=1}^{m}\left(1-\frac{\sqrt{b_j^2 +\epsilon_j^2}}{\sqrt{|\va_j^*\vz|^2 +\epsilon_j^2}}\right)\va_j\va_j^*\vz.
	\end{align*}
	The magnitude of $\nabla f_{\vepsilon}(\vz)$ is under control even when $ |\va_j^*\vz|/|\va_j^*\vx| $ very small. This fact avoids the extreme value of gradients during each update, which makes each update flows in a desirable direction and guarantees the gradient satisfies \emph{curvature condition}. The curvature condition shall be introduced in Lemma~\ref{curvature condition}.
	
	So the truncation-based methods (TAF, TWF) use truncation to withdraw the spurious components and RAF uses re-weighted to reduce the effects of ``bad'' gradients. Compared to them, our PAF controls these components by adding the perturbed term, i.e., $\epsilon$ to avoid the extreme value during each update, which frees our methods from truncation or re-weighted procedure. Besides, such a perturbation and corresponding benefit is applicable to both real and complex fields, thus make our theoretical analysis easily incorporate the complex field as a whole.
	
    Numerical tests show that our proposed algorithm outperforms AF ($ \vepsilon=\0 $) in terms of success rate for real signals, as shown in Figure \ref{success_res}.
	Besides, using vanilla gradient descent to solve the perturbed amplitude-based model (\ref{model}), we can achieve linear convergence with  $ m=\mathcal{O}(n) $ measurements for both real and complex signals  (see Section II).  The result improves upon the WF method, which uses $ m=\mathcal{O}(n\log n) $ measurements, or the AF method, which can be theoretically proved only for real signals, or the TWF, TAF, RAF methods, which need truncation or re-weighted procedure during each iteration.

	In summary, compared with the previous algorithms for solving model (\ref{amplitude_model}) or (\ref{likelihood_model}), the PAF method needn't truncation or re-weighted at all and the convergence result holds for both real and complex signals. Numerical experiments show that the proposed PAF method is slightly more efficient although comparable computationally with TAF, RAF and significantly more efficient than TWF (see Section III). We believe the reason lies in the fact that truncated/re-weighted methods, such as TWF, TAF, RAF incur additional computational cost on measuring the gradient components.
	
	\subsection{Notations}
	Let $ \vx\in\H^n $ $ (\H=\C $ or $ \H =\R) $  be the target  signal. Throughout this paper, we assume that $ \va_j \in\H^n$, $ j=1,\ldots,m $ are $ m $ independent and identically distributed standard Gaussian random measurement vectors, i.e. $ \va_j \sim \mathcal{N}(0,I) $ for $ \H=\R $ and $ \va_j \sim \mathcal{N}(0,I/2) + i\mathcal{N}(0,I/2)  $ for $ \H=\C $. For each measurement $ \va_j $, we obtain $ b_j = |\va_j^*\vx| $.  We shall attempt to recover the original signal $ \vx $ from $ b_j $, $ j=1,\ldots,m $ by solving the perturbed amplitude-based model (\ref{model}). In  this paper, we use $ C $, $ c $ or the subscript/superscript form of them to represent constants and their values vary according to the context. Since for phase retrieval the best we can do is to recover the target signal $ \vx $ up to a global phase/sign, we use the following definition for distance between two vectors $\vx,\vz\in\H^n$:
	\begin{equation}\label{distance}
	{\rm dist} (\vz,\vx)=\min_{\phi\in[0,2\pi)}\|\vz - e^{i\phi}\vx\|:=\|\vz - e^{i\phi_\vx(\vz)}\vx\|,
	\end{equation}
	where
	\begin{equation}  \label{phi_argmin}
	\phi_\vx (\vz) := \argmin{\phi\in[0,2\pi)} \|\vz - e^{i\phi}\vx\|.
	\end{equation}
	For any $\rho\geq 0$, we define the $ \rho $-neighborhood of $ \vx $ as
	\begin{equation}\label{neighbour}
	\mathcal{S}_\vx(\rho):=\Big\{\vz\in\C^d: \text{dist}(\vz,\vx)\leq\rho\|\vx\|\Bigr\}.
	\end{equation}
	
	% =========================Sec===============================

	\section{Perturbed Amplitude Flow Algorithm}
	
	\subsection{Initialization}\label{initial_part}
	To avoid iterations getting trapped in undesirable stationary points, a proper initialization is essential to any non-convex optimization problem. To achieve this goal, many initialization methods have been proposed, such as the spectral initialization method \cite{candes_wf}, a modified spectral initialization method \cite{chen2015solving} and the null initialization method \cite{wang2018solving}. These methods are all based on finding the eigenvector corresponding to the largest eigenvalue of a specially designed Hermitian matrix.
	
	Here we adopt the initialization strategy given in \cite{gao2017phaseless}, which is shown to provide a good initial guess under $ \mathcal{O}(n) $ measurements. With this strategy, the initial guess $ \vz_0 $ is obtained by calculating the eigenvector corresponding to the largest eigenvalue of the Hermitian matrix
	$$
	Y = \frac{1}{m}\sum_{j=1}^{m}\left( \gamma - \exp(-b^2_j/\lambda^2)\right)\va_j\va_j^*
	$$
	with $ \gamma = 1/2$ for $\H = \C$ or $ \gamma = 1/\sqrt{3}$ for $\H = \R$, and normalized to $ \|\vz_0\|=\lambda $, where $ \lambda $ is defined by
	$$
	\lambda^2 = \frac{1}{m}\sum_{j=1}^{m}b^2_j.
	$$
	\begin{lemma}[\hspace{-1sp}\cite{gao2017phaseless}]\label{init}
		Let $\vz_0$ be the above initial guess. For any $ \xi>0 $, there exists a $C_\xi>0$ such that for $ m\geq C_\xi n $,
		$$
		\textup{dist}(\vz_0, \vx) \leq \xi\|\vx\|
		$$
		holds with probability at least $ 1-4\exp(-c_\xi n) $.
	\end{lemma}
	
	\subsection{Gradient Descent Iteration}\label{gra_section}
	
	After initialization to obtain $\vz_0$, we use gradient descent on the loss function $f_{\vepsilon}$ given in (\ref{model}) by
	$$
	f_{\vepsilon}(\vz):= \frac{1}{m}\sum_{j=1}^{m}\left(\sqrt{|\va_j^*\vz|^2+\epsilon_j^2}
	- \sqrt{ b_j^2+\epsilon_j^2} \right)^2
	$$
	to iteratively refine the estimation:
	\begin{equation}\label{iteration}
	\vz_{k+1} = \vz_{k} - \mu \nabla f_{\vepsilon}(\vz_{k}),
	\end{equation}
	where $\mu$ is the step size and  $\nabla f_{\vepsilon}(\vz) $ is the Wirtinger derivative of $ f_{\vepsilon}(\vz) $ with respect to $ \vz $ in complex variables $\vz,\overline{\vz}$ which is defined as
	\begin{align*}
	\nabla f_{\vepsilon}(\vz) &:= \left(\frac{\partial f_{\vepsilon}(\vz,\overline{\vz})}{\partial \vz}\Big|_{\overline{\vz}
		= \text{constant}}\right)^* \\
	& = \frac{1}{m}\sum_{j=1}^{m}\left(1-\frac{\sqrt{b_j^2 +\epsilon_j^2}}{\sqrt{|\va_j^*\vz|^2 +\epsilon_j^2}}\right)\va_j\va_j^*\vz.
	\end{align*}
	As simple as the scheme  (\ref{iteration}) may look, our main result proves that it can achieve linear convergence under the optimal order of measurements $ m = \mathcal{O}(n) $ by choosing $ \vepsilon=\sqrt{\alpha} \vb $ for an appropriately chosen parameter $ \alpha>0 $ ($0.37 \leq \alpha \leq 29$).
	
	Motivated by the technique used in WF, the proof of our main result is mainly based on the following two key lemmas, whose proofs are given in Section \ref{lemma_proof_sec}.
	\begin{lemma}\label{smoothness}
		Let $\vx$ be the target signal and assume that $ \vepsilon$ satisfies (\ref{condition_eps}). For any $ \delta>0$, there exist constants $C_\delta$, $c_\delta >0$ such that as long as $m \geq C_\delta n$, then with probability at least $ 1-\exp(-c_\delta n)$,
		\begin{equation}\label{smooth}
		\|\nabla f_{\vepsilon}(\vz)\|\leq (1+\delta)\cdot\textup{dist}(\vz, \vx)
		\end{equation}
		holds for every $z \in \mathbb{H}^n$ satisfying $z \in \mathcal{S}_\vx(1/10)$.
	\end{lemma}
	
	This lemma implies that the gradient of $f_{\vepsilon}$ is well controlled in the neighborhood of the target signal $\vx$.
	
	\begin{lemma}\label{curvature condition}
		Let $\vx$ be the target signal and assume that $ \vepsilon = \sqrt{\alpha}\vb $ with  $ 0.37\leq \alpha\leq 29 $. There exist positive constants $C,c, \beta_\alpha$ depending on $\alpha$ such that for any $ \vz\in \mathcal{S}_\vx(1/10) $ and $m\geq Cn$, we have
		\begin{equation}\label{curvature}
		\Re\big(\langle \nabla f_{\vepsilon}(\vz),\, \vz-\vx e^{i\phi_\vx(\vz)}\rangle\big)
		\geq \beta_\alpha\cdot\textup{dist}^2(\vz, \vx)
		\end{equation}
		with probability at least $ 1-\exp(-c n) $.
	\end{lemma}
		
	The constants in the lemma can, in theory, be explicitly estimated, although the theoretical estimates are typically ``overkills'' for practical applications, just like in other existing schemes. Later in Remark \ref{relation},
	we show more explicitly the relation between $\beta_\alpha$ and $\alpha$. Particularly,  by setting $  \alpha =0.826 $,  $ \beta_\alpha = 64/5945 $ roughly reaches its largest value.     For $ \vepsilon = \sqrt{\alpha}\vb  $ with $ \alpha\in[0.37, 29] $,  Lemma \ref{curvature condition} guarantees sufficient descent along the search direction.
	
	Set $ \vh := e^{-i\phi_\vx(\vz)}\vz-\vx $ with $ \rho=\|\vh\|$.  Then
	\begin{small}
		\begin{align*}
		&\Re\big(\langle \nabla f_{\vepsilon}(\vz),\, \vz-\vx e^{i\phi_\vx(\vz)}\rangle\big) \\
		&=\frac{1}{m}\sum_{j=1}^{m} \left(1-\frac{\sqrt{b_j^2+\epsilon_j^2}}{\sqrt{|\va_j^*(\vx+\vh)|^2+\epsilon_j^2}}\right)\big(|\va_j^*\vh|^2+\Re(\vh^*\va_j\va_j^*\vx) \big).	
		\end{align*}
	\end{small}	
	The main technique in proving Lemma \ref{curvature condition} is that we first fix one $ \vz\in\C^n $ and then provide estimates separately for cases  $ |\va_j^*\vh| \geq \rho |\va_j^*\vx|$ and $|\va_j^*\vh|< \rho|\va_j^*\vx| $. An $\eta  $-net argument is then used to obtain uniform control over all $ \vz\in\SS_\vx(\rho) $.
	
	Building on these two lemmas, we can now state and prove our main theorem, which establishes linear convergence of the PAF algorithm iteration (\ref{iteration}).
	
	\begin{theorem}\label{linear_converge}
		Under the conditions of Lemma~\ref{curvature condition},  let $\vz_{k}$, $k\in {\mathbb Z}_{+}$ be the iterations generated by (\ref{iteration}) with $ \mu=
\beta_\alpha/1.001^2 $. Assume that $\vz_0 \in \mathcal{S}_\vx(1/10) $. Then there
exist positive constants $ C, c $ such that for $ m\geq Cn $, with probability at
least $ 1-\exp(-c n) $, the following holds for all $k\in {\mathbb Z}_{+}$
		\[
		\textup{dist}^2(\vz_{k+1}, \vx)\leq (1-\beta_\alpha^2/1.001^2)\cdot \textup{dist}^2(\vz_k, \vx).
		\]
		In particular by taking $ \alpha=0.826 $, with probability at least $
1-\exp(-c n) $, the following holds for all $k\in {\mathbb Z}_{+}$
		\[
		\textup{dist}(\vz_{k}, \vx)\leq \frac{1}{10}
		\left(1-\frac{0.0107^2}{1.001^2}\right)^{k/2}\cdot \|\vx\|.
		\]
	\end{theorem}
	\begin{IEEEproof}
		According to the update rule (\ref{iteration}), Lemma \ref{smoothness} and Lemma \ref{curvature condition}, for $ m\geq C n $, with probability at least $ 1-\exp(-cn) $ we have
		\begin{align*}
		&\textup{dist}^2(\vz_{k+1}, \vx)\\
		&\leq\|\vz_{k+1} - \vx e^{i\phi_\vx(\vz_k)}\|^2\\
		&=\| \vz_k - \vx e^{i\phi_\vx(\vz_k)}-\mu \nabla f_{\vepsilon}(\vz_k)\|^2\\
		& =\|\vz_{k}-\vx e^{i\phi_\vx(\vz_k)}\|^2 \\
		&\quad- 2\mu\Re\big(\langle \nabla f_{\vepsilon}(\vz_k), \,\vz_k-\vx e^{i\phi_\vx(\vz_k)}\rangle\big)+\mu^2 \|\nabla f_{\vepsilon}(\vz_k)\|^2\\
		&\leq \|\vz_{k}-\vx e^{i\phi_\vx(\vz_k)}\|^2 - 2\mu\cdot \beta_\alpha\|\vz_{k}-\vx e^{i\phi_\vx(\vz_k)}\|^2\\
		&\quad \quad +\mu^2\cdot 1.001^2\|\vz_{k}-\vx e^{i\phi_\vx(\vz_k)}\|^2\\
		&=\big(1-\mu \cdot(2\beta_\alpha- 1.001^2\mu)\big)\|\vz_{k}-\vx e^{i\phi_\vx(\vz_k)}\|^2\\
		&= (1-\beta_\alpha^2/1.001^2)\cdot \textup{dist}^2(\vz_k, \vx).
		\end{align*}
		This establishes the linear convergence part of the theorem.
		
		For the second part, we set $\alpha=0.826$. Later in Remark \ref{relation}, we show that one may take $ \beta_\alpha=64/5945 $ in $ \mu= \beta_\alpha/1.001^2 $. Substituting these values in we thus obtain 		
		\begin{align*}
		\textup{dist}(\vz_{k}, \vx)&\leq (1-\beta_\alpha^2/1.001^2)^{1/2}\cdot \textup{dist}(\vz_{k-1}, \vx)\\
		&< (1-0.0107^2/1.001^2)^{1/2}\cdot\textup{dist}(\vz_{k-1}, \vx)\\
		&\leq\frac{1}{10} \left(1-\frac{0.0107^2}{1.001^2}\right)^{k/2}\cdot\| \vx\|.
		\end{align*}
	\end{IEEEproof}
	
	As mentioned earlier, we can achieve $\vz_0 \in \mathcal{S}_\vx(1/10)$ through initialization given in Lemma \ref{init}, by setting $\xi = 1/10$. This also requires $ m=\mathcal{O}(n)$ measurements. Thus the combination of Lemma \ref{init} and Theorem \ref{linear_converge} yield linear convergence of the PAF algorithm.

	%================================================
	
	\section{Numerical Experiments} 	
	\subsection{Simulation Study}
		
	To evaluate the performance of our PAF algorithm, we present a series of simulated tests and compare them with WF, TWF, AF, TAF and RAF. We perform all the simulations under the same initialization procedure. All experiments are carried out on Matlab 2017b with a 2.3 GHz Intel Core i5-8259U and 16 GB memory.

\begin{figure}[h]
	\graphicspath{{figures/}}
	\begin{center}
		\subfigure[]{
			\includegraphics[width=0.48\textwidth]{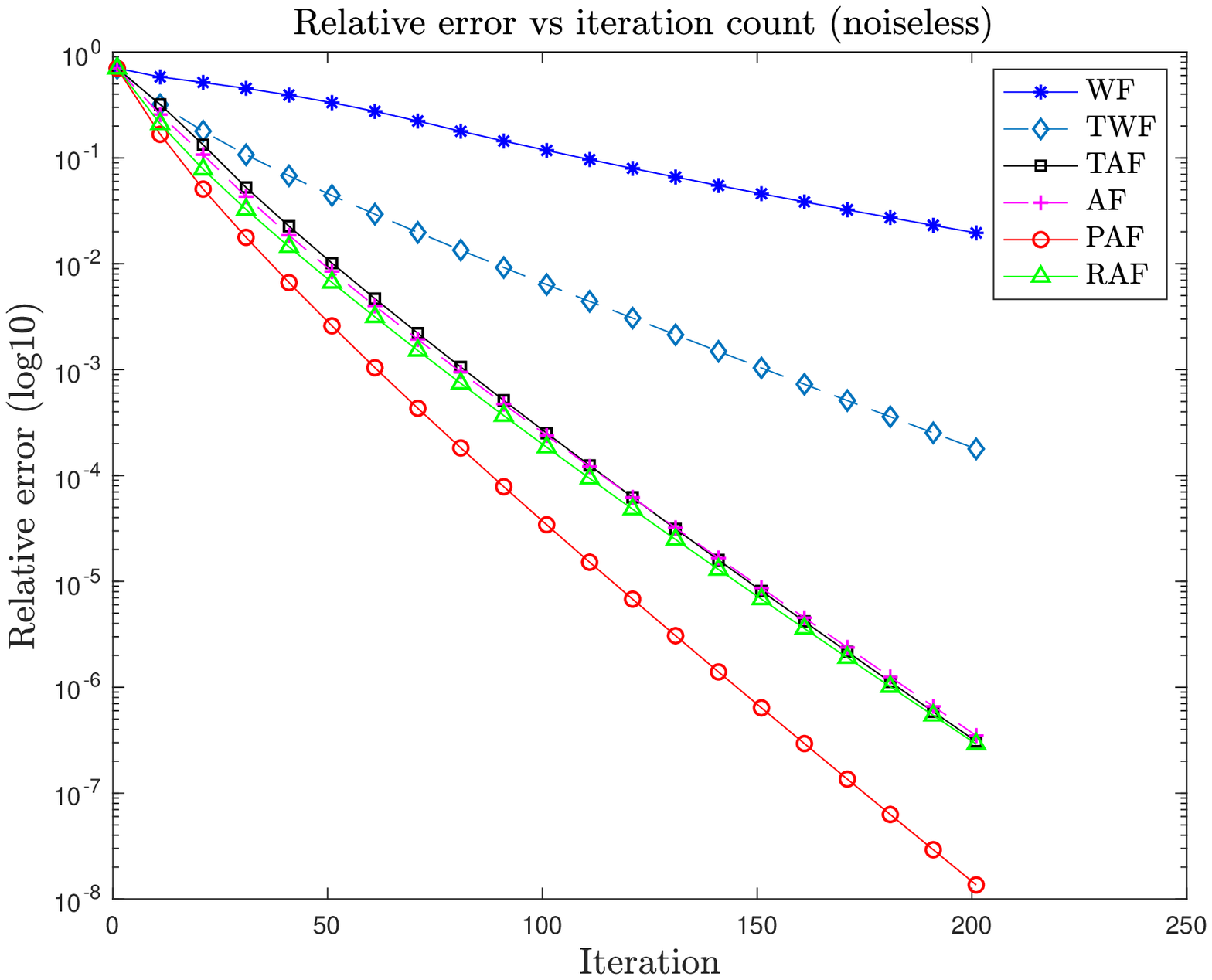}}
		\subfigure[]{
			\includegraphics[width=0.48\textwidth]{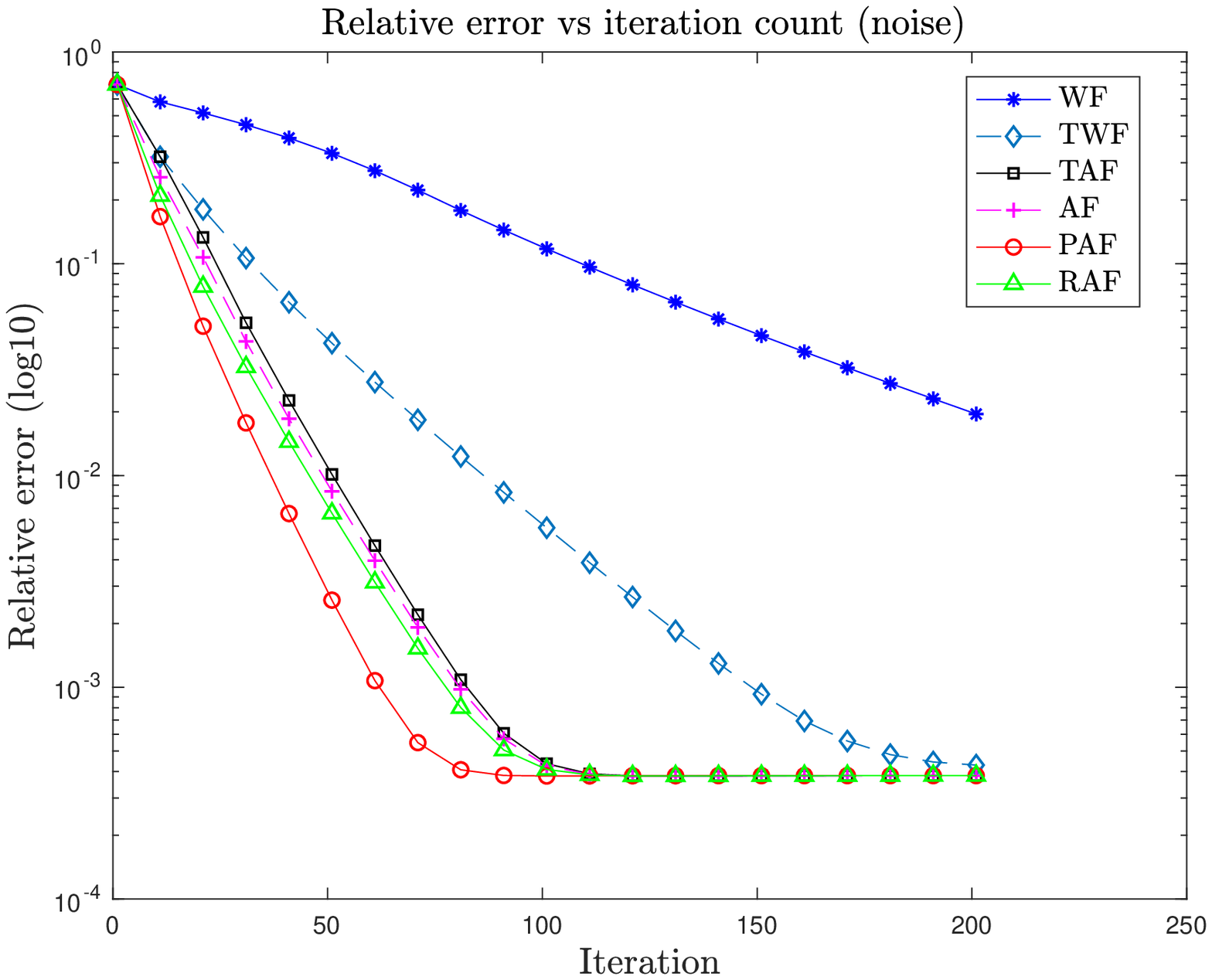}}
		\caption{\em
			Convergence experiments: Plot of relative error (log(10)) vs number of iterations  for PAF (our algorithm), WF, TWF, TAF and AF method. Take $n = 512, m = 4.5n$.  The figure (a) (for the exact measurements) and figure (b) (for noisy measurements ) both show that PAF method provides better
			solution and also converges faster. }\label{compare_algoris}
	\end{center}
\end{figure}

First we plot the relative error for the recovery of a complex-valued signal, in logarithmic scale versus the iteration count for WF, TWF, AF, TAF, RAF and PAF. We choose $ n = 512 $ with $ m = 4.5n $ i.i.d. Gaussian random measurements $ \va_1,\va_2,\ldots,\va_m\in\C^n$. For the initialization, we follow the method given in Section \ref{initial_part} with 50 power iterations. For the PAF algorithm we set $ \vepsilon=\vb $ and fix the step size $ \mu = 2.5$. Note that AF is equivalent to PAF algorithm with $ \vepsilon=\0 $. We also consider the case where the measurements are contaminated by noise, i.e. $ \vb=\abs{A\vx}+\omega $ where the noise $ \omega$ follows distribution $ \omega\sim \mathcal{N}(0,I/10) $.
The results are plotted in Figure \ref{compare_algoris}. It shows that PAF, TWF, AF, TAF and RAF, all of which converge linearly in theory, have comparable convergence rate. PAF seems to have a slight advantage possibly due to its ability to handle a larger step size.

\begin{figure}[h]
	\graphicspath{{figures/}}
	\begin{center}
		\subfigure[]{
			\includegraphics[width=0.48\textwidth]{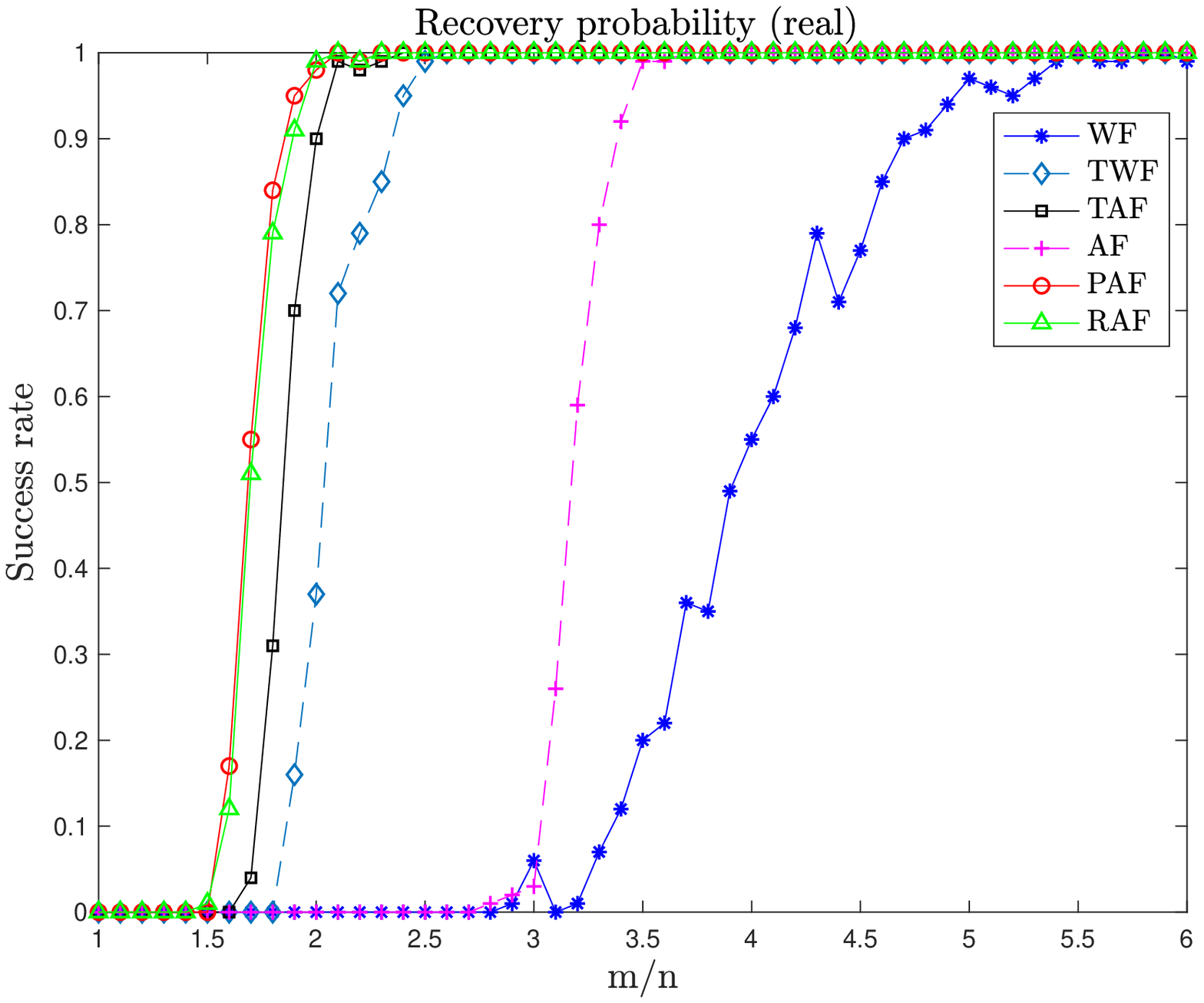}}
		\subfigure[]{
			\includegraphics[width=0.48\textwidth]{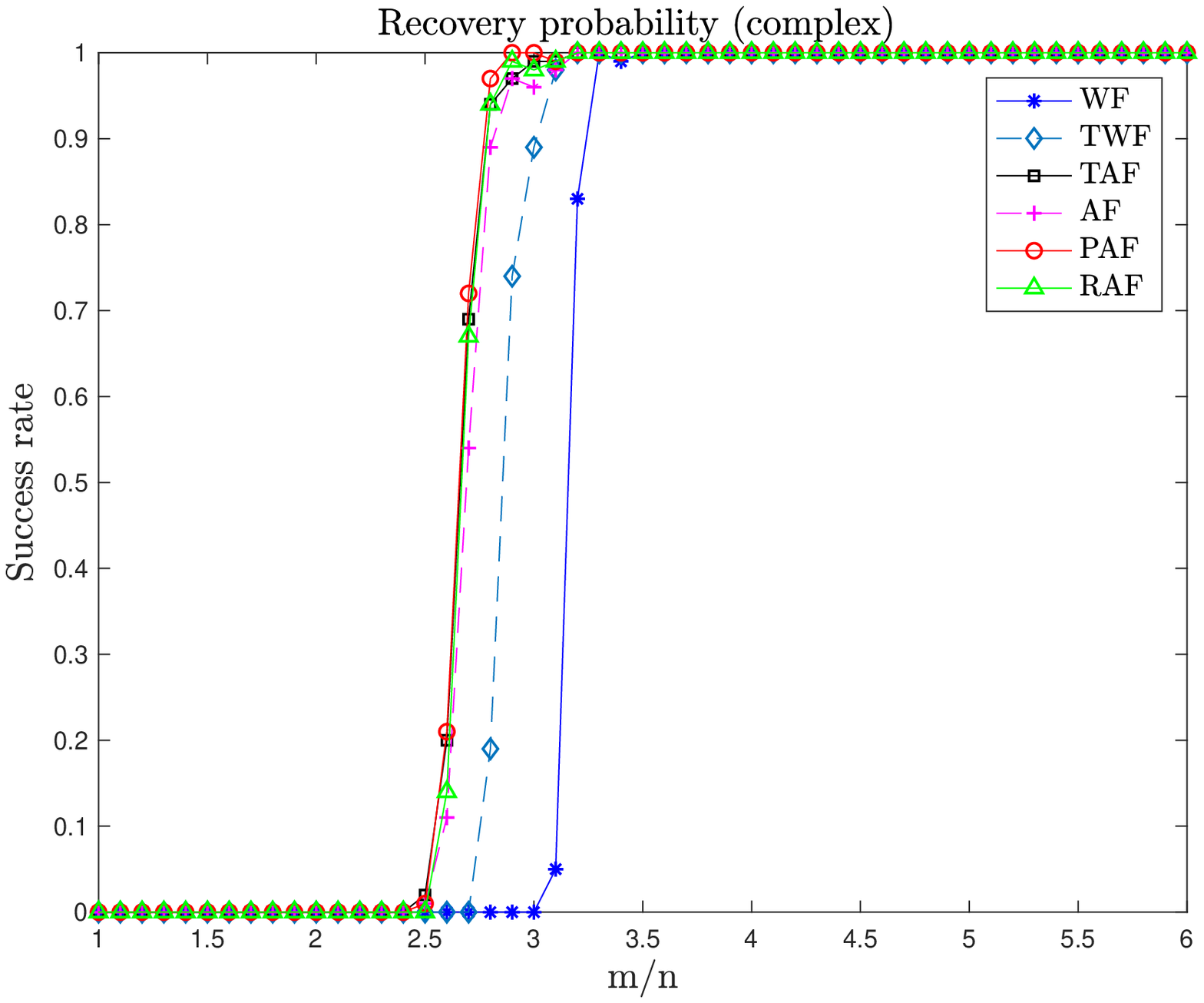}}
		\caption{ \em Success rate experiments: Empirical probability of successful test based on 100 random trails for different $m/n$. Take $n = 512$ and change $m/n$ between 1 and 6.
			The figures demonstrate that PAF, TWF, TAF and AF are better than WF in terms of the success rate.
		}\label{success_res}
	\end{center}
\end{figure}

Next, we compare the empirical success rate of PAF with that of
WF, TWF, AF, TAF and RAF. Here we set the maximum number of gradient-type iterations to $ T= 2500 $ for each scheme.  In PAF, we set $n=512$, $ \vepsilon=\vb $ and fix the step size to $ \mu = 1 $.
We let $ m/n $ vary  from $1$ to $6$. A test is successful if the relative error is within $ 10^{-5} $ after the maximum number of iterations.  For the test we compute the success rate by performing 100 random trials for each $ m/n $. The results are given in  Figure \ref{success_res}. Of particular note is that in the real case, PAF, TWF and TAF all perform better than AF, indicating the effectiveness of controlling the size of the gradient in all gradient descent algorithms for avoiding spurious stationary points. WF seems to lag behind other algorithms, unsurprisingly, as it agrees with the theoretical analysis.

\subsection{Recovery of Natural Image}

To show the efficiency and scalability of our algorithm, we use PAF to recover the Milky Way Galaxy image \footnote{Download from http://pics-about-space.com/milky-way-galaxy}, which is the image used in \cite{candes_wf,cdp} with the coded diffraction measurements. We denote the image by $\bm{X}$,  $ \bm{X}\in\R^{1080\times1920\times 3} $. This is a color image so it has three channels. Thus we actually perform phase retrieval for each of the three channels separately. Let $\vx$ denote any of the color channels of $\bm{X}$. We have measurements
$$
\vb^{(l)} = \abs{\bm{F}\bm{D}^{(l)}\vx},\,\, 1\leq l\leq L,
$$
where $ \bm{F} $ denotes the $ n\times n $ discrete Fourier transform matrix, and $ \bm{D}^{(l)} $ is a diagonal matrix having i.i.d. entries sampled from a distribution $ g $. Here we take the {\it octanary } pattern that $ g=g_1g_2$, where $ g_1 $ and $ g_2 $ are independent with distributions
$$
g_1=\left\{
\begin{aligned}
1\quad & \text{with prob.} & 1/4 \\
-1\quad & \text{with prob.} & 1/4 \\
-i\quad & \text{with prob.} & 1/4 \\
i\quad & \text{with prob.} & 1/4
\end{aligned}
\right.
$$
and
$$
g_2=\left\{
\begin{aligned}
\sqrt{2}/2 \quad & \text{with prob.} & 4/5 \\
\sqrt{3}\quad & \text{with prob.} & 1/5
\end{aligned}
\right..
$$
We set $L=20$ and adopt the same initialization method for all schemes in our comparison. For each model, we record the time elapsed and the iterations needed to achieve relative error at $ 10^{-5} $ and $ 10^{-10} $, respectively. The results are shown in Table \ref{table1}. It is shown that PAF achieves the same level of precision and is comparable in efficiency with AF and TAF. Besides, note that it took TAF, RAF, PAF and AF the same number of iterations to achieve fixed relative error. Moreover, it's reasonable that our PAF is a little bit slower than AF ($ \vepsilon=\0 $) with additional nonzero item $ \vepsilon $.
These three methods are significantly more efficient than WF and TWF.
\begin{table}[h]	
	\centering
	\fontsize{10}{13}\selectfont
	\begin{threeparttable}
		\caption{Iteraions and elapsed time.}
		\label{table1}
		\begin{tabular}{cccc}
			\toprule
			Algorithm & Relative error & Iter & Time(s) \cr
			\cmidrule(lr){1-4}
			\multirow{2}{*}{WF}& $ 1\times10^{-5} $ & 172 & 348.36\cr
			&  $ 1\times10^{-10} $ & 302 &  606.71  \cr
			\cmidrule(lr){1-4}
			\multirow{2}{*}{TWF}&  $ 1\times10^{-5} $& 51 & 320.53\cr
			& $ 1\times10^{-10} $ & 118 & 691.39 \cr
			\cmidrule(lr){1-4}
			\multirow{2}{*}{TAF}&  $ 1\times10^{-5} $& 37 & 118.06\cr
			& $ 1\times10^{-10} $& 84 & 250.23 \cr
			\cmidrule(lr){1-4}
			\multirow{2}{*}{RAF}&  $ 1\times10^{-5} $&  37  &  124.99\cr
			& $ 1\times10^{-10} $ & 84  & 271.40 \cr
			\cmidrule(lr){1-4}
			\multirow{2}{*}{PAF}&  $ 1\times10^{-5} $&  37 & 97.27 \cr
			& $ 1\times10^{-10} $ & 84 & 224.23 \cr
			\cmidrule(lr){1-4}
			\multirow{2}{*}{AF}&  $ 1\times10^{-5} $&  37 & 87.65 \cr
			& $ 1\times10^{-10} $ & 84 & 189.24 \cr
			\bottomrule
		\end{tabular}
	\end{threeparttable}
\end{table}
\begin{figure}
	\graphicspath{{figures/}}
	\centering
	\includegraphics[width=3in,height=1.2in]{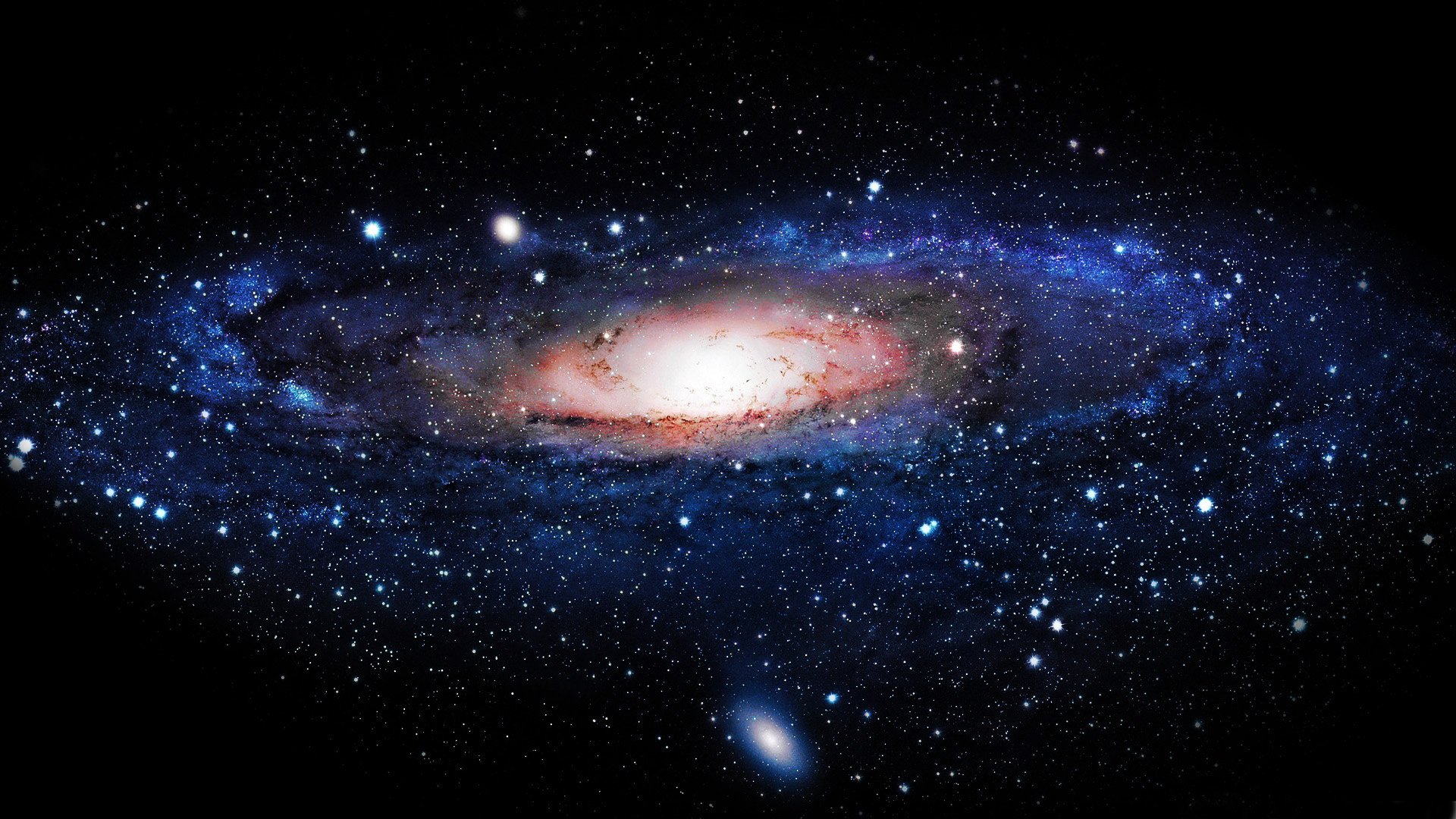}
	\caption{\em Milky Way Galaxy image: Image size is $ 1080\times 1920$ pixels. Our PAF algorithm with $ L=6 $ takes  300
		%200
		iterations, computation time is 183.5 sec, %156
		%	331.7 sec,
		relative error is
		%$ 3.1\times10^{-16} $.
		$ 5.04\times10^{-15}  $. }
	\label{fig:galaxy}
\end{figure}

Interestingly if we take a much smaller $L=6$, while WF does not recover the target image, our PAF method actually performs better than with $L=20$. It takes  300 iterations and computation time $ 183.5 $ sec to achieve recovery with a relative error of $ 5.04\times10^{-15}$ in Figure \ref{fig:galaxy}.
While more iterations are taken here, the computational time is actually less because $L=6$ is significantly smaller than $L=20$.
	
	%==================================
	
	\section{Proof of main lemmas in section \ref{gra_section}}\label{lemma_proof_sec}

	\subsection{Proof of Lemma \ref{smoothness}}
	
	\begin{IEEEproof}
		For any $ \vz\in\C^n $, set $ \vh = e^{-i\phi_\vx(\vz)} \vz -\vx $, where we recall that $\phi_\vx(\vz)$ is given in (\ref{phi_argmin}). Then $\|\vh\|=\textup{dist}(\vz,\vx)$.
		Denote $ A = [\va_1, \ldots, \va_m]^*\in\C^{m\times n}$,  $ \vv=[v_1, v_2,\ldots, v_m]^T$
		with $ v_j =\left(1-\frac{\sqrt{b_j^2 +\epsilon_j^2}}{\sqrt{|\va_j^*\vz|^2 +\epsilon_j^2}}\right)(\va_j^*\vz)  $. Note that we set $v_j=0$ if $b_j=\epsilon_j = \va_j^*\vz=0$. Then $\nabla f_{\vepsilon}(\vz)=\frac{1}{m}A^*\vv$. For any $ \epsilon_j>0 $, we have
		\begin{align*}
		|v_j|^2 &= \left|1-\frac{\sqrt{b_j^2 +\epsilon_j^2}}{\sqrt{|\va_j^*\vz|^2 +\epsilon_j^2}}\right|^2|\va_j^*\vz|^2\\
		&=\frac{\left(\sqrt{|\va_j^*\vz|^2 +\epsilon_j^2}-\sqrt{|\va_j^*\vx|^2 +\epsilon_j^2}\right)^2}{|\va_j^*\vz|^2 +\epsilon_j^2}|\va_j^*\vz|^2\\
		&\leq \left(\sqrt{|\va_j^*\vz|^2 +\epsilon_j^2}-\sqrt{|\va_j^*\vx|^2 +\epsilon_j^2}\right)^2\\
		&=\left(\sqrt{|\va_j^*(\vx+\vh\label{\vz})|^2 +\epsilon_j^2}-\sqrt{|\va_j^*\vx|^2 +\epsilon_j^2}\right)^2\\
		&\leq |\va_j^*\vh|^2,
		\end{align*}
		where the last inequality follows from the inequality $|\sqrt{t^2+c^2}-\sqrt{s^2+c^2}| \leq |t-s|$ for any $t,s,c\in\R$. According to Lemma \ref{sub_gaussian_concentration} (see the Appendix), for any $ \delta'>0 $ and $ m\geq C_{\delta'}n $ with a sufficiently large constant $ C_{\delta'}$, the inequality
		$$
		\|\vv\|^2 =\sum_{j=1}^{m}\left|v_j\right|^2\leq\sum_{j=1}^{m}\left|\va_j^*\vh\right|^2\leq(1+\delta')m\|\vh\|^2
		$$
		holds with probability at least $ 1-e^{- c_{\delta'}n} $ for some $c_{\delta'}>0$. Also for the Gaussian random matrix $ A $ and any $\delta''>0$, for $ m\geq C_{\delta''}n $ we have $ \|A^*\|\leq(1+\delta'')\sqrt{m}  $ with probability at least $ 1-e^{-c_{\delta''}n} $ (\cite{vershynin2010introduction}, Remark 5.40). These results together imply that
		\begin{align*}
		\|\nabla f_{\vepsilon}(\vz)\|&=\frac{1}{m}\|A^*\vv\|\\
		&\leq\frac{1}{m}\|A^*\|\|\vv\|\\
		&\leq\sqrt{(1+\delta')}(1+\delta'')\|\vh\|\\
		&\leq(1+\delta)\|\vh\|
		\end{align*}
		holds with probability at least $ 1-\exp(-c_\delta n) $ whenever $ m\geq C_\delta n $ for some $C_\delta, c_\delta>0$. Here we choose $ 1+\delta\geq\sqrt{(1+\delta')}(1+\delta'')$ and $ C_\delta\geq\max\{ C_{\delta'}, C_{\delta''}\} $.
	\end{IEEEproof}
	
	\subsection{Proof of Lemma \ref{curvature condition}}
	
	\begin{IEEEproof}
		Without loss of generality, we shall assume that the target signal $\vx$ has $ \|\vx\|=1 $. Again for each $ \vz\in\C^n $ we set $ \vh = e^{-i\phi_\vx(\vz)} \vz -\vx $, and denote $ \tilde{\vh} = \vh/\|\vh\| $.  Definition \ref{distance} implies that $ \Im (\vh^*\vx)=0 $. Since $ \vz\in \mathcal{S}_\vx(1/10)  $,  we have $ \rho:=\|\vh\|\leq 1/10 $. Therefore
		\begingroup
		\allowdisplaybreaks
		\begin{small}
			\begin{align*}
			&\Re\big(\langle \nabla f_{\vepsilon}(\vz),\, \vz-\vx e^{i\phi_\vx(\vz)}\rangle\big)\\
			&=\Re\big(\langle \nabla f_{\vepsilon}(\vz), \, e^{i\phi_\vx(\vz)}\vh \rangle\big)	\\
			&= \frac{1}{m}\sum_{j=1}^{m}\left(1-\frac{\sqrt{b_j^2 +\epsilon_j^2}}{\sqrt{|\va_j^*\vz|^2 +\epsilon_j^2}}\right)\Re\big((\va_j^*\vz)\,e^{-i\phi_\vx(\vz)}(\vh^*\va_j)\big)\\
			&=\frac{1}{m}\sum_{j=1}^{m}\frac{\sqrt{|\va_j^*\vz|^2 +\epsilon_j^2}-\sqrt{|\va_j^*\vx|^2 +\epsilon_j^2}}{\sqrt{|\va_j^*\vz|^2 +\epsilon_j^2}}\Re\big(\va_j^*(\vx+\vh)(\vh^*\va_j)\big)\\
			&=\frac{1}{m}\sum_{j=1}^{m}\frac{\big(|\va_j^*\vz|^2 -|\va_j^*\vx|^2\big)\Re\big(\va_j^*(\vx+\vh)(\vh^*\va_j)\big) }{\sqrt{|\va_j^*\vz|^2 +\epsilon_j^2}\left(\sqrt{|\va_j^*\vz|^2 +\epsilon_j^2}+\sqrt{|\va_j^*\vx|^2 +\epsilon_j^2}\right)}\\
			&=\frac{1}{m}\sum_{j=1}^{m}\frac{2\big(\Re(\vh^*\va_j\va_j^*\vx)\big)^2+ 3\Re(\vh^*\va_j\va_j^*\vx)|\va_j^*\vh|^2+|\va_j^*\vh|^4}{\sqrt{|\va_j^*\vz|^2 +\epsilon_j^2}\left(\sqrt{|\va_j^*\vz|^2 +\epsilon_j^2}+\sqrt{|\va_j^*\vx|^2 +\epsilon_j^2}\right)}\\
			&=\frac{1}{m}\sum_{j=1}^{m}T_j,
			\end{align*}
		\end{small}	
	    \endgroup
		with $ T_j $ being the $ j $-th item of the summation. To simplify the statement, we use $ d_j $ to denote the denominator of $ T_j $, i.e.,
		$$
		d_j =\sqrt{|\va_j^*\vz|^2 +\epsilon_j^2}\left(\sqrt{|\va_j^*\vz|^2 +\epsilon_j^2}+\sqrt{|\va_j^*\vx|^2 +\epsilon_j^2}\right).
		$$
		To prove the conclusion holds for all $ \vz\in\mathcal{S}_\vx(1/10) $, i.e., any $ \tilde{\vh}  $ in unit ball. We first consider $ \tilde{\vh}\in \C^n $ to be fixed and then divide it into two cases.
		
		In the first case, we assume $\tilde{\vh}=c\vx$ with $\abs{c}=1$. Here we have $\Im(\tilde{\vh}^*\vx)=0$, which implies $ \tilde{\vh}= \pm \vx$. Hence
		\begingroup
		\allowdisplaybreaks
		\begin{align*}
		d_j&=\sqrt{|\va_j^*\vz|^2 +\epsilon_j^2}\left(\sqrt{|\va_j^*\vz|^2 +\epsilon_j^2}+\sqrt{|\va_j^*\vx|^2 +\epsilon_j^2}\right)\\
		&\leq\left(3(1+\|\vh\|^2)+2\alpha+1/2\right)|\va_j^*\vx|^2\\
		&\leq (353/100+2\alpha)|\va_j^*\vx|^2,
		\end{align*}
		\endgroup
		due to the facts that
		$$
		|\va_j^*\vz|^2=|\va_j^*(\vx+\vh)|^2\leq 2\big(|\va_j^*\vx|^2+\|\vh\|^2\cdot|\va_j^*\vx|^2\big),
		$$
		$ \epsilon_j^2=\alpha |\va_j^*\vx|^2 $ and  $ a+\sqrt{ab}\leq \frac{3}{2}a+\frac{1}{2}b $.
		Thus under the condition of  $ \|\vh\|\leq\frac{1}{10} $, we obtain
		\begingroup
		\allowdisplaybreaks
		\begin{align*}
		T_j
		&=\frac{\big(2\pm 3\|\vh\|+\|\vh\|^2\big)|\va_j^*\vx|^4}{d_j}\|\vh\|^2\\
		&\geq \frac{\big(2\pm 3\|\vh\|+\|\vh\|^2\big)|\va_j^*\vx|^2}{353/100+2\alpha}\|\vh\|^2\\
		&\geq \frac{171}{353+200\alpha}|\va_j^*\vx|^2\|\vh\|^2.
		\end{align*}
		\endgroup
		By Lemma \ref{sub_gaussian_concentration} of the Appendix, for $ m\geq C_\delta n $, with probability greater than $ 1-\exp(-c_\delta m) $ we have
%		\begin{equation}
%		\begin{aligned}\label{specialcase}
%		& \Re\Big(\langle \nabla f_{\vepsilon}(\vz),\, \vz-\vx e^{i\phi(\vz)}\rangle\Big)  =\frac{1}{m}\sum_{j=1}^{m} T_j \\
%		&\geq \frac{1}{m}\sum_{j=1}^{m} \frac{171}{353+200\alpha}|\va_j^*\vx|^2\|\vh\|^2\\
%		&\geq  \frac{171}{353+200\alpha}(1-\delta)\|\vh\|^2.
%		\end{aligned}
%		\end{equation}
	    \begin{equation}
		\begin{aligned}\label{specialcase}
		& \Re\Big(\langle \nabla f_{\vepsilon}(\vz),\, \vz-\vx e^{i\phi(\vz)}\rangle\Big) \\
		& =\frac{1}{m}\sum_{j=1}^{m} T_j \\
		&\geq \frac{1}{m}\sum_{j=1}^{m} \frac{171}{353+200\alpha}|\va_j^*\vx|^2\|\vh\|^2\\
		&\geq  \frac{171}{353+200\alpha}(1-\delta)\|\vh\|^2.
		\end{aligned}
		\end{equation}
		
		For the second case $ \tilde{\vh}\neq \pm \vx $, given the assumption $\|\vx\|=1$ and $\|\vh\|=\rho$, we claim that
		\begin{equation}  \label{pro3}
		\PP(\rho|\va_j^*\vx|> |\va_j^*\vh|)=\PP(\rho|\va_j^*\vx|\leq |\va_j^*\vh|)=1/2.
		\end{equation}
		Indeed, for each measurement $ \va_j $ we have
		\begin{equation}\label{pro1}
		\PP(\rho|\va_j^*\vx|= |\va_j^*\vh|)=\PP(|\va_j^*\vx|= |\va_j^*\tilde{\vh}|)=0.
		\end{equation}
		Also note that a Gaussian random measurement $\va$ is rotational invariant, i.e. for any unitary matrix $ O $, $O\va $ is also a Gaussian random measurement. Thus for fixed $ \vx $ and $\tilde{\vh} $, we may without loss of generality assume that $\tilde{\vh}=\ve_1 $ and $ \vx=\sigma \ve_1 + \sqrt{1-\sigma^2}\ve_2$, with $ \sigma=\tilde{\vh}^*\vx\in\R $. This is because otherwise we can always find a unitary matrix to map $\tilde{\vh}, \vx$ to these two vectors. Set
		$$
		O:=\begin{pmatrix}
		O_1 & \0\\
		\0& O_2\\
		\end{pmatrix}\in\C^{n\times n}
		$$
		where $ O_2 \in\C^{(n-2)\times (n-2)}$ is unitary and	
		$$
		O_1 = \begin{pmatrix}
		\sigma & \sqrt{1-\sigma^2}\\
		\sqrt{1-\sigma^2}&-\sigma\\
		\end{pmatrix}.
		$$
		Then we have $O\vx = \tilde{\vh}$ and $ O\tilde{\vh}=\vx $.  Set $ \vg := O\va $  and $\vg$ is a Gaussian random measurement. Consequently we have
		\begin{align*}
		\PP(|\va^*\vx|>|\va^*\tilde{\vh}|)& = \PP(|\vg^*O\vx|>|\vg^*O\tilde{\vh}|) \\
		&= \PP(|\vg^*\tilde{\vh}|>|\vg^*\vx|),
		\end{align*}
		which  implies
		\begin{equation} \label{pro2}
		\PP(\rho|\va_j^*\vx|> |\va_j^*\vh|)=\PP(\rho|\va_j^*\vx|< |\va_j^*\vh|).
		\end{equation}
		Combining  (\ref{pro1}) and (\ref{pro2}) we now obtain (\ref{pro3}).

		For each index set $ I\subseteq \{1,2,\ldots,m\} $, define a corresponding event
		$$
		\E_I:=\bigl\{ \rho\,|\va_j^*\vx|>|\va_j^*\vh|, \,\,\forall j\in I ;\,\, \rho\,|\va_k^*\vx|\leq|\va_k^*\vh|, \,\,\forall k\in I^c\bigr\} .
		$$
		According to (\ref{pro3}), we know that the event $ \E_I  $ occurs with probability $ 1/2^m $.  We assume that  $ I_0 $ is an index set which satisfies $ \frac{m}{4}\leq|I_0|\leq \frac{3m}{4} $.  Then on event $ \E_{I_0} $, $\Re\big(\langle \nabla f_{\vepsilon}(\vz),\, \vz-\vx e^{i\phi_\vx(\vz)}\rangle\big) $
		can be divided into two  groups:
		\begin{align*}
		m\,\Re\big(\langle \nabla f_{\vepsilon}(\vz), \, \vz-e^{i\phi_\vx(\vz)}\vx \rangle\big)
		=\sum_{j\in I_0}T_j+\sum_{k\in I_0^c}T_k.
		\end{align*}
		For each group, we next provide an upper bound and a lower bound for the denominators $ d_j $, $ j=1,\ldots,m$. Recall that $ \vepsilon = \sqrt{\alpha}\vb \,(\alpha>0) $. When $ j\in I_0 = \big\{j\,:\, \rho\,|\va_j^*\vx|>|\va_j^*\vh|\big\} $  we have
		\begin{equation}
			\begin{aligned} \label{1upperbound}
			d_j & = \sqrt{|\va_j^*\vz|^2 +\epsilon_j^2}\left(\sqrt{|\va_j^*\vz|^2 +\epsilon_j^2}
			+\sqrt{|\va_j^*\vx|^2 +\epsilon_j^2}\right) \\  & \leq \frac{3}{2}|\va_j^*\vz|^2 + 2\epsilon_j^2
			+ \frac{1}{2}|\va_j^*\vx|^2    \\
			 & < \frac{3}{2}(1+\rho)^2|\va_j^*\vx|^2 + 2\alpha|\va_j^*\vx|^2
			+\frac{1}{2}|\va_j^*\vx|^2  \\
			& = \left(2\alpha + 2 +3\rho+\frac{3}{2}\rho^2\right)|\va_j^*\vx|^2 \\
	        & = U_1|\va_j^*\vx|^2
			\end{aligned}
		\end{equation}
		where $ U_1:= 2\alpha + 2 + 3\rho+\frac{3}{2}\rho^2 $. Here the second inequality follows from $ \epsilon_j^2 = \alpha|\va_j^*\vx|^2 $ and
		\begin{align*}
		|\va_j^*\vz|^2&=|\va_j^*(\vx+\vh)|^2\\
		&\leq(\abs{\va_j^*\vx}+\abs{\va_j^*\vh})^2\\
		&<(1+\rho)^2\abs{\va_j^*\vx}^2 .
		\end{align*}
		On the other hand, since
		\begin{align*}
		|\va_j^*\vz|^2&=|\va_j^*(\vx+\vh)|^2\\
		&\geq(\abs{\va_j^*\vx}-\abs{\va_j^*\vh})^2\\
		&>(1/\rho-1)^2\abs{\va_j^*\vh}^2
		\end{align*}
		and $ \epsilon_j^2=\alpha\abs{\va_j^*\vx}^2>(\alpha/\rho^2)\abs{\va_j^*\vh}^2 $, we have
		\begin{small}
		\begin{equation}\label{1lowerbound}
			\begin{aligned}
			d_j &= \sqrt{|\va_j^*\vz|^2 +\epsilon_j^2}\left(\sqrt{|\va_j^*\vz|^2 +\epsilon_j^2}+\sqrt{|\va_j^*\vx|^2 +\epsilon_j^2}\right)\\
			&> \sqrt{(1/\rho-1)^2+(\alpha/\rho^2)} \Big( \sqrt{(1/\rho-1)^2+(\alpha/\rho^2)}\\
			&\quad\cdot|\va_j^*\vh|+(\sqrt{1+\alpha}/\rho)\cdot|\va_j^*\vh|\Big)|\va_j^*\vh|\\		
			&=\frac{\sqrt{(1-\rho)^2+\alpha}\,\big(\sqrt{(1-\rho)^2+\alpha}+\sqrt{1+\alpha}\big)}{\rho^2}
			|\va_j^*\vh|^2\\
			&= L_1|\va_j^*\vh|^2,
			\end{aligned}
		\end{equation}
		\end{small}
		where $ L_1:=  \sqrt{(1-\rho)^2+\alpha}\,\big(\sqrt{(1-\rho)^2+\alpha}+\sqrt{1+\alpha}\big)/\rho^2$.		
		Similarly, for $ k\in  I_0^c = \big\{k\,:\, \rho\,|\va_k^*\vx|\leq|\va_k^*\vh|\big\} $, we have
		$$
		|\va_k^*\vz|^2\leq\big(|\va_k^*\vx|+|\va_k^*\vh|\big)^2\leq(1+1/\rho)^2 |\va_k^*\vh|^2,
		$$
		and hence
		\begin{equation}
			\begin{aligned}\label{2upperbound}
			d_k &=\sqrt{|\va_k^*\vz|^2 +\epsilon_k^2}\left(\sqrt{|\va_k^*\vz|^2 +\epsilon_k^2}+\sqrt{|\va_k^*\vx|^2 +\epsilon_k^2}\right)\\
			&\leq \frac{3}{2}|\va_k^*\vz|^2 + 2\epsilon_k^2 + \frac{1}{2}|\va_k^*\vx|^2   \\
			&\leq\frac{3}{2}(1/\rho+1)^2|\va_k^*\vh|^2+(2\alpha+1/2)/\rho^2\cdot|\va_k^*\vh|^2\\
			&=\left(\frac{2\alpha+2}{\rho^2}+\frac{3}{2}+\frac{3}{\rho}\right)|\va_k^*\vh|^2\\
			&= U_2|\va_k^*\vh|^2,
			\end{aligned}
		\end{equation}
		where $ U_2 :=\frac{2\alpha+2}{\rho^2}+\frac{3}{2}+\frac{3}{\rho}$ and
		\begin{equation}
		\begingroup
        \allowdisplaybreaks
			\begin{aligned}\label{2lowerbound}
			d_k &=\sqrt{|\va_k^*\vz|^2 +\epsilon_k^2}\left(\sqrt{|\va_k^*\vz|^2 +\epsilon_k^2}+\sqrt{|\va_k^*\vx|^2 +\epsilon_k^2}\right)\\
			&\geq\epsilon_k\Big(\epsilon_k+\sqrt{|\va_k^*\vx|^2 +\epsilon_k^2}  \Big)\\
			&=\left(\alpha+\sqrt{\alpha(1+\alpha)}  \right)|\va_k^*\vx|^2\\
			&=L_2|\va_k^*\vx|^2,
			\end{aligned}
		\endgroup
		\end{equation}			
		where $ L_2 := \alpha+\sqrt{\alpha(1+\alpha)} $.
		
		Using the concentration inequalities given in the Appendix, we next give the  lower bounds of $ \sum_{j\in I_0}T_j $ and $ \sum_{j\in I_0^c}T_k $. Based on (\ref{1upperbound}),
(\ref{1lowerbound}) and Lemma \ref{concentration}, given any $ \delta>0 $, for $
|I_0| \geq C_1(\delta)n $ the following inequality holds with probability at least $
1-\exp\big(-c_1(\delta)\cdot|I_0|\big) $
		\begin{small}
		\begingroup
        \allowdisplaybreaks
			\begin{align*}
			& \sum_{j\in I_0}T_j  \\
			&= \sum_{j\in I_0}\Bigg(\frac{\big(\sqrt{2}\Re(\vh^*\va_j\va_j^*\vx)+\frac{3}{2\sqrt{2}}|\va_j^*\vh|^2\big)^2}{d_j} - \frac{|\va_j^*\vh|^4}{8d_j}\Bigg)\\
			&\geq\sum_{j\in I_0} \Bigg(\frac{\big(\sqrt{2}\Re(\vh^*\va_j\va_j^*\vx)+\frac{3}{2\sqrt{2}}|\va_j^*\vh|^2\big)^2}{U_1|\va_j^*\vx|^2} - \frac{|\va_j^*\vh|^4}{8L_1|\va_j^*\vh|^2}\Bigg)\\
			%			& =\sum_{j\in I_0} \Bigg(\frac{2\big(\Re(\vh^*\va_j\va_j^*\vx)\big)^2+ 3\Re(\vh^*\va_j\va_j^*\vx)|\va_j^*\vh|^2+\frac{9}{8}|\va_j^*\vh|^4}{U_1|\va_j^*\vx|^2} \\
			%			&\quad \quad - \frac{|\va_j^*\vh|^2}{8L_1}\Bigg)\\
			& \geq \sum_{j\in I_0} \Bigg(\frac{2\big(\Re(\vh^*\va_j\va_j^*\vx)\big)^2- 3|\vh^*\va_j\va_j^*\vx||\va_j^*\vh|^2}{U_1|\va_j^*\vx|^2} - \frac{|\va_j^*\vh|^2}{8L_1}\Bigg)\\
			&\geq\sum_{j\in I_0}\Bigg(\frac{2}{U_1}\frac{\big(\Re(\vh^*\va_j\va_j^*\vx)\big)^2}{|\va_j^*\vx|^2} - \frac{3\|\vh\|^3}{U_1}|\va_j^*\tilde{\vh}|^2-\frac{|\va_j^*\vh|^2}{8L_1}\Bigg)\\
			 &\geq|I_0|\cdot\|\vh\|^2\Bigg(\frac{2}{U_1}\Big(\frac{1}{8}+\frac{7}{32}\Re^2(\tilde{\vh}^*\vx)\Big)-\frac{3}{2U_1}\|\vh\|-\frac{1}{16L_1}-\frac{\delta}{4}\Bigg)\\
			%			 &\geq|I_0|\cdot\|\vh\|^2\left(\frac{1}{4U_1}-\frac{3\rho}{2U_1}-\frac{1}{16L_1}+\frac{7}{16U_1}\Re^2(\tilde{\vh}^*\vx)-\frac{\delta}{4}\right)\\
			 &\geq|I_0|\cdot\|\vh\|^2\left(\frac{1}{4U_1}-\frac{3\rho}{2U_1}-\frac{1}{16L_1}-\frac{\delta}{4}\right)\\
			&=|I_0|\cdot\|\vh\|^2\cdot \varphi_1,
			\end{align*}
		\endgroup
		\end{small}	
		where $ \varphi_1:=\frac{1-6\rho}{4U_1}-\frac{1}{16L_1}-\frac{\delta}{4} $.
Here the fourth inequality comes from   Lemma \ref{concentration}.
		
		Similarly, according to (\ref{2upperbound}), (\ref{2lowerbound}) and Lemma
\ref{concentration},  for $ |I_0^c|\geq C_2(\delta)n $ the following inequality holds
 with probability at least $ 1-\exp\big(-c_2(\delta)\cdot|I_0^c|\big) $:
 \begin{small}
		\begingroup
        \allowdisplaybreaks
			\begin{align*}
			&\sum_{k\in I_0^c} T_k \\
			&= \sum_{k\in I_0^c} \Bigg( \frac{\big(\frac{3}{2}\Re(\vh^*\va_k\va_k^*\vx)+|\va_k^*\vh|^2\big)^2}{d_k}-\frac{\big(\Re(\vh^*\va_k\va_k^*\vx)\big)^2}{4d_k} \Bigg)\\
			&\geq \sum_{k\in I_0^c}\Bigg(\frac{\frac{9}{4}\big(\Re(\vh^*\va_k\va_k^*\vx)\big)^2+ 3\Re(\vh^*\va_k\va_k^*\vx)|\va_k^*\vh|^2+|\va_k^*\vh|^4}{U_2|\va_k^*\vh|^2}\\
			&\quad\quad-\frac{\big(\Re(\vh^*\va_k\va_k^*\vx)\big)^2}{4L_2|\va_k^*\vx|^2}\Bigg )\\
			&\geq|I_0^c|\cdot \Bigg(\frac{9}{4U_2}\Big(\frac{1}{8}+\frac{7}{32}\Re^2(\tilde{\vh}^*\vx)\Big)+\frac{3\|\vh\|}{2U_2}\Re(\tilde{\vh}^*\vx)\\
			 &\quad\quad+\frac{\|\vh\|^2}{2U_2}-\frac{\|\vh\|^2}{4L_2}\Big(\frac{3}{8}+\frac{9}{32}\Re^2(\tilde{\vh}^*\vx)\Big)-\frac{\delta\|\vh\|^2}{4}\Bigg)\\
			&\geq|I_0^c|\cdot \|\vh\|^2\Bigg(\frac{9}{32U_2\rho^2}+\frac{1}{2U_2}-\frac{3}{32L_2}\\
			 &\quad+\Big(\frac{63}{128U_2\rho^2}-\frac{9}{128L_2}\Big)\Re^2(\tilde{\vh}^*\vx)+\frac{3}{2U_2\rho}\Re(\tilde{\vh}^*\vx)-\frac{\delta}{4} \Bigg)\\
			&\geq |I_0^c|\cdot \|\vh\|^2\Bigg(\frac{9}{32U_2\rho^2}+\frac{1}{2U_2}-\frac{3}{32L_2}-\phi-\frac{\delta}{4} \Bigg)\\
			&=|I_0^c|\cdot \|\vh\|^2\cdot \varphi_2,
			\end{align*}
		\endgroup
		\end{small}	
		where $ \phi = \Big(\frac{3}{4U_2\rho}\Big)^2/
\Big(\frac{63}{128U_2\rho^2}-\frac{9}{128L_2}\Big)$ and
		$
\varphi_2:=\frac{9}{32U_2\rho^2}+\frac{1}{2U_2}-\frac{3}{32L_2}-\phi-\frac{\delta}{4}$.
The second inequality follows from the concentration inequalities given in Lemma
\ref{concentration}. The fourth inequality derives from the facts that $
\frac{63}{128U_2\rho^2}-\frac{9}{128L_2}>0$ for any $ 0.37\leq\alpha\leq 197$ and $
\rho \leq 1/10 $.
		
		Set $ \delta:=0.001 $.	For arbitrary fixed $ \alpha\in[0.37, 197] $, a simple observation is that $ \varphi_1 $ and $ \varphi_2 $ are decreasing functions of  $ \rho $.  So we next only consider $ \rho = 1/10 $. When $ 0.37\leq\alpha\leq 197$, we have
		\begin{equation}\label{r1}
		\varphi_1=\frac{1}{10U_1}-\frac{1}{16L_1}-\frac{\delta}{4}>0
		\end{equation}
		and
		\begin{equation}\label{r2}
		\varphi_2=\frac{229}{8U_2}-\frac{3}{32L_2}+\frac{225}{4U_2^2\tilde{\phi}}-\frac{\delta}{4} >0,
		\end{equation}
		with $ \tilde{\phi} = \frac{9}{128L_2}-\frac{1575}{32U_2}<0$.
		
		For sufficiently large constant $ C\geq 4\max\{C_1(\delta), C_2(\delta)\} $, as long as $ m\geq Cn $, we have $ |I_0|\geq m/4 \geq C_1(\delta)n $ and $ |I_0^c|\geq m/4 \geq C_2(\delta)n $. Thus with probability at least $ (1-\exp(-c_3m))/2^m $,  we have
		\begin{equation}
			\begin{aligned}\label{generalcase}
			&\Re\big(\langle \nabla f_{\vepsilon}(\vz),\, \vz-\vx e^{i\phi_\vx(\vz)}\rangle\big)\\
			&=\frac{1}{m}\Big(\sum_{j\in I_0}T_j+\sum_{k\in I_0^c}T_k\Big)\\
			&\geq\frac{1}{m}\bigg(|I_0|\|\vh\|^2\cdot \varphi_1+|I_0^c|\|\vh\|^2\cdot \varphi_2 \bigg)\\
			&\geq\frac{1}{m}\left(\frac{m}{4}\|\vh\|^2\cdot \varphi_1+\frac{m}{4}\|\vh\|^2\cdot \varphi_2\right)\\
			&=(\varphi_1+\varphi_2)\|\vh\|^2/4.
			\end{aligned}
		\end{equation}
		The number of the index sets $I$  satisfying $  \frac{m}{4}\leq |I|\leq \frac{3m}{4}$ is $ \sum_{k=m/4}^{3m/4} {m\choose k}$.  So for fixed $ \tilde{\vh} $, when $ \tilde{\vh}\neq\pm \vx $, the inequality (\ref{generalcase})
		holds with probability greater than $ \sum_{k=m/4}^{3m/4} {m\choose k}(1-\exp(-c_3m))/2^m $. Note that
		\begin{align*}
		&\sum_{3m/4+1}^m{m\choose k}=\sum_{k=0}^{m/4-1}{m\choose k}=\sum_{k=0}^{m/4-1}{m\choose m/4-1-k}\\
		&={m\choose m/4-1}\sum_{k=0}^{m/4-1}\frac{(m/4-1)\cdots (m/4-k)}{(3m/4+2)\cdots (3m/4+k+1)}\\
		&<{m\choose m/4-1}\sum_{k=0}^{m/4-1} \left(\frac{m/4-1}{3m/4+2}\right)^k\\
		&< (4e)^{m/4} \cdot \frac{3}{2},
		\end{align*}
		and $(4e)^{1/4}<2$. Hence $\sum_{k=0}^{m/4-1}{m\choose k}/2^m<c_0^m$ for some $c_0\in (0,1)$, which implies that $ \sum_{k=m/4}^{3m/4} {m\choose k}(1-\exp(-c_3m))/2^m\geq 1-\exp(-c_5m) $. Moreover, for $  \alpha \in [0.37, 197] $ we have
		\begin{equation}\label{comparison}
		\frac{\varphi_1+\varphi_2}{4} < \frac{171}{353+200\alpha}\cdot(1-0.001).
		\end{equation}        	           	
		Considering the two cases as a whole, for a fixed $ \vz $, combining (\ref{specialcase}), (\ref{generalcase}) and (\ref{comparison}), we obtain
		\begin{equation}\label{aimfunc}
		\Re\big(\langle \nabla f_{\vepsilon}(\vz),\, \vz-\vx e^{i\phi_\vx(\vz)}\rangle\big)
		\geq \frac{\varphi_1+\varphi_2}{4}\|\vh\|^2
		\end{equation}
		with probability at least $ 1-\exp(-c_6m)$. Particularly, when $ \alpha\in[0.37, 29] $ we have $  \frac{\varphi_1+\varphi_2}{4}>0.001$.
		
		To complete the proof, we will need to establish uniform bound over all vectors, so we adopt an $\eta $-net argument. Observe that
		\begin{align*}
		&\Re\big(\langle \nabla f_{\vepsilon}(\vz),\, \vz-\vx e^{i\phi_\vx(\vz)}\rangle\big)\\
		&=\Re\big(\langle \nabla f_{\vepsilon}(e^{-i\phi_\vx(\vz)}\vz),\, e^{-i\phi_\vx(\vz)}\vz-\vx\rangle\big)\\
		&= \Re\big(\langle \nabla f_\vepsilon (\vx+\rho\tilde{\vh}),\, \rho\tilde{\vh} \rangle\big).
		\end{align*}
		For any $ \vz\in\C^n $, which means for any $ \tilde{\vh} $ with $ \|\tilde{\vh}\|=1 $ and $ \Im(\tilde{\vh}^*\vx)=0 $, we consider the function $ \Re\big(\langle \nabla f_\vepsilon (\vx+\rho\tilde{\vh}),\, \rho\tilde{\vh} \rangle\big) $ with $ \rho\leq 1/10 $. Suppose that $ \tilde{\vh}_1, \tilde{\vh}_2\in \C^n $ satisfy  $ \|\tilde{\vh}_1-\tilde{\vh}_2\|\leq \eta $. When  $ 0.37 \leq\alpha\leq 29 $ we have
		\begin{align*}
		&\big|\Re\big(\langle \nabla f_\vepsilon (\vx+\rho\tilde{\vh}_1), \rho\tilde{\vh}_1 \rangle\big)-\Re\big(\langle \nabla f_\vepsilon (\vx+\rho\tilde{\vh}_2), \rho\tilde{\vh}_2 \rangle\big)\big| \\
		& \leq\big|\Re\big(\langle \nabla f_\vepsilon (\vx+\rho\tilde{\vh}_1), \rho(\tilde{\vh}_1-\tilde{\vh}_2)\rangle\big)\big|\\
		&\quad +\big|\Re\big(\langle \nabla f_\vepsilon (\vx+\rho\tilde{\vh}_1) - \nabla f_\vepsilon (\vx+\rho\tilde{\vh}_2), \rho\tilde{\vh}_2 \rangle\big)\big|\\
		&\leq \rho\|\nabla f_\vepsilon (\vx+\rho\tilde{\vh}_1)\|\cdot \eta + \|\nabla^2 f_\vepsilon (\xi)\|\cdot\rho^2\|\tilde{\vh}_2\|\cdot \eta\\
		&\leq 2\rho^2\|\tilde{\vh}_1\|\cdot \eta + 2\sqrt{\frac{1+\alpha}{\alpha}}\cdot\rho^2\|\tilde{\vh}_2\|\cdot \eta\\
		&=2\Bigg(1+\sqrt{\frac{1+\alpha}{\alpha}}\Bigg) \cdot \eta \cdot \rho^2\\
		&<6\eta \cdot \rho^2,
		\end{align*}
		where $ \xi\in \C^n $. Here the third inequality follows from Lemma \ref{smoothness} and Lemma \ref{hessian}. Therefore for any $ \tilde{\vh}_1 $ and $ \tilde{\vh}_2 $ satisfying $ \|\tilde{\vh}_1 - \tilde{\vh}_2\|\leq \eta:=\frac{\delta}{6} $ with $ \delta = 0.001 $, let $ \NN_\eta $ be an $ \eta $-net for the unit sphere of $ \C^n $ with cardinality $ |\NN_\eta|\leq (1+2/\eta)^{2n}$. Then for all $ \vz $, $ 0.37\leq\alpha\leq 29 $ and $m\geq (C_2\cdot \eta^{-2}\log \eta^{-1})n$, with probability at least $ 1-\exp(-c n) $ we have
		\begin{align*}
		\Re\big(\langle \nabla f_{\vepsilon}(\vz),\, \vz-\vx e^{i\phi_\vx(\vz)}\rangle\big)&\geq \big((\varphi_1+\varphi_2)/4-\delta\big)\|\vh\|^2\\
		&=\beta_\alpha\|\vh\|^2
		\end{align*}
		with $ \beta_\alpha:= (\varphi_1+\varphi_2)/4-\delta >0$.  According to Remark \ref{relation}, when $ \alpha=0.826$, $ \beta_\alpha=64/5945 $ approximately reaches its largest value. 	
	\end{IEEEproof}
	
	\begin{remark}\label{relation}
		According to the proof of Lemma \ref{curvature condition},
		by taking $\rho=1/10$ and $\delta=0.001$,  we have
		$U_1=2\alpha+463/200$, $U_2=200\alpha + 463/2$, $L_1=100\alpha + 81+100\sqrt{(\alpha+1)(\alpha+0.81)}$ and $L_2=\alpha+\sqrt{\alpha(1+\alpha)}$.
		Recall that
		\begin{align*}
		\beta_\alpha &= (\varphi_1+\varphi_2)/4-\delta\\
		 &=\frac{1}{40U_1}+\frac{229}{32U_2}-\frac{1}{64L_1}-\frac{3}{128L_2}+\frac{225}{16U_2^2\tilde{\phi}}-\frac{9\delta}{8},
		\end{align*}
		with $ \tilde{\phi} = \frac{9}{128L_2}-\frac{1575}{32U_2}$. Figure \ref{re} here shows the relationship between $ \beta_\alpha $ and  $ \alpha $.
		\begin{figure}[H]
			\graphicspath{{figures/}}
			\begin{center}
				{\includegraphics[width=0.48\textwidth]{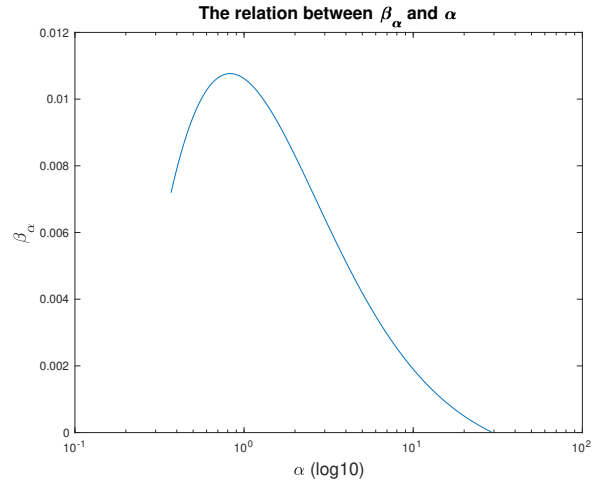}}
			\end{center}
			\caption{The relationship between $ \beta_\alpha $ and $ \alpha (\log 10) $.}\label{re}
		\end{figure}
		Particularly,  when $ \alpha=0.826$, $ \beta_\alpha=64/5945$  roughly reaches its maximum.
	\end{remark}
		
	%-----------------------------

	\appendix[Auxiliary Lemmas]
	
	In previous sections we have applied concentration inequalities several times. They have played a key role in the proof of our results. Here we present these concentration inequalities used for the proof of Lemma \ref{smoothness} and Lemma \ref{curvature condition}.

	\begin{lemma}[\cite{candes2013phaselift} Lemma 3.1 ]\label{sub_gaussian_concentration}
		Let $ \va_1,\va_2,\ldots,\va_m \in\C^n$ be i.i.d. Gaussian random measurements. Fix any $ \delta $ in $ (0,1/2) $ and assume $ m\geq 20\delta^{-2}n $. Then for all unit vectors $ \vu\in\C^n $,
		\[
		1-\delta\leq\frac{1}{m}\sum_{j=1}^{m}|\va_j^*\vu|^2\leq 1+\delta
		\]
		holds with probability at least $ 1-\exp(-m t^2/2) $, where $ \delta/4 = t^2+t $.
	\end{lemma}
	\begin{lemma}\label{expectation}
		Let $ \va\in\C^n$  be a Gaussian random measurement. Let $ \vx\in\C^n $ and $ \tilde{\vh}\in\C^n $ be two fixed vectors with $ \|\vx\|=\|\tilde{\vh}\|=1 $, $ \Im(\tilde{\vh}^*\vx)=0 $ and $ \tilde{\vh}\neq\pm \vx $. Then we have
		\begin{equation}
			\begin{aligned}\label{exp-haax}
			&\E\big(\Re(\tilde{\vh}^*\va\va^*\vx)\cdot I_{\{|\va^*\vx|>|\va^*\tilde{\vh}|\}} \big)\\
			&=\E\big(\Re(\tilde{\vh}^*\va\va^*\vx)\cdot I_{\{|\va^*\vx|\leq|\va^*\tilde{\vh}|\}} \big)\\
			&=\Re(\tilde{\vh}^*\vx),
			\end{aligned}
		\end{equation}		
		\begin{align}\label{exp-ax1}
		\frac{1}{2}\leq\E\big(|\va^*\vx|^2\cdot I_{\{|\va^*\vx|>|\va^*\tilde{\vh}|\}}\big)\leq\frac{3}{4},
		\end{align}
		\begin{align}\label{exp-ax2}
		\frac{1}{4}\leq\E\big(|\va^*\vx|^2\cdot I_{\{|\va^*\vx|\leq|\va^*\tilde{\vh}|\}}\big)\leq\frac{1}{2},
		\end{align}
%		\begin{small}
		\begin{equation}
		\begin{aligned}\label{exp-haax/ax1}
		\frac{1}{8}+\frac{7}{32}\Re^2(\tilde{\vh}^*\vx)&\leq\E\bigg(\frac{\big(\Re(\tilde{\vh}^*\va\va^*\vx)\big)^2}{|\va^*\vx|^2}\cdot I_{\{|\va^*\vx|>|
			\va^*\tilde{\vh}\}}\bigg)\\
		&\leq\frac{1}{4}+\frac{1}{4}\Re^2(\tilde{\vh}^*\vx)
		\end{aligned}
		\end{equation}
%		\end{small}
		and
		\begin{equation}
		\begin{aligned}\label{exp-haax/ax2}
		 \frac{1}{4}+\frac{1}{4}\Re^2(\tilde{\vh}^*\vx)&\leq\E\bigg(\frac{\big(\Re(\tilde{\vh}^*\va\va^*\vx)\big)^2}{|\va^*\vx|^2}\cdot I_{\{|\va^*\vx|\leq|
			\va^*\tilde{\vh}\}}\bigg)\\
		&\leq\frac{3}{8}+\frac{9}{32}\Re^2(\tilde{\vh}^*\vx).
		\end{aligned}
		\end{equation}
	\end{lemma}
	\begin{IEEEproof}
		Since the distribution of $ \va $ is invariant by unitary transformation, we can take $ \vx=\ve_1 $ and $ \tilde{\vh}=\sigma\ve_1+\sqrt{1-\sigma^2}\ve_2 $, where $ \sigma = \vx^*\tilde{\vh}=\Re( \vx^*\tilde{\vh})\in\R $ and $ |\sigma|<1 $. We use $ \xi_1, \xi_2, \xi_3, \xi_4$ to represent the real and imaginary parts of $ a_1 $ and $ a_2 $ respectively, which implies that the variables $ \xi_1, \xi_2, \xi_3, \xi_4 $ are independent and obey normal distribution $ \mathcal{N}(0,1/2) $.  Then it follows that
		\begin{align*}
		 \E\big(\Re(\tilde{\vh}^*\va\va^*\vx)\big)&=\E\big(\sigma(\xi_1^2+\xi_2^2)+\sqrt{1-\sigma^2}(\xi_1\xi_3+\xi_2\xi_4)\big)\\
		&=\sigma=\Re(\tilde{\vh}^*\vx),
		\end{align*}
		\[
		\E(|\va^*\vx|^2)=\E(\xi_1^2+\xi_2^2)=1
		\]
		and
		\begin{align*}
		&\E\bigg(\frac{\big(\Re(\tilde{\vh}^*\va\va^*\vx)\big)^2}{|\va^*\vx|^2}\bigg)\\
		 &=\E\bigg(\frac{\big(\sigma(\xi_1^2+\xi_2^2)+\sqrt{1-\sigma^2}(\xi_1\xi_3+\xi_2\xi_4)\big)^2}{\xi_1^2+\xi_2^2}\bigg)\\
		&=\frac{1}{2}+\frac{1}{2}\sigma^2.
		\end{align*}
		Since  $ \va $ is invariant by unitary transformation and $ \vx, \tilde{\vh}$ are two fixed vectors satisfying $ \tilde{\vh}\neq \pm\vx $, so we have
		\begin{align*}
		&\E\big(\Re(\tilde{\vh}^*\va\va^*\vx)\cdot I_{\{|\va^*\vx|\leq|\va^*\tilde{\vh}|\}}\big)\\
		&=\E\big(\Re(\tilde{\vh}^*\va\va^*\vx)\cdot I_{\{|\va^*\vx|<|\va^*\tilde{\vh}|\}}\big)\\
		&=\E\big(\Re(\vx^*\vg\vg^*\tilde{\vh})\cdot I_{\{|\vg^*\vx|>|\vg^*\tilde{\vh}|\}}\big)\\
		&=\E\big(\Re(\tilde{\vh}^*\va\va^*\vx)\cdot I_{\{|\va^*\vx|>|\va^*\tilde{\vh}|\}}\big).
		\end{align*}
		Here $ \vg:=O\va $ is a Gaussian random measurement with unitray matrix $ O $ satisfying $ O\vx = \tilde{\vh} $ and $ O\tilde{\vh} = \vx $. Then we obtain
		\begin{align*}
		& \E\big(\Re(\tilde{\vh}^*\va\va^*\vx)\big)=\E\big(\Re(\tilde{\vh}^*\va\va^*\vx)\cdot I_{\{|\va^*\vx|>|\va^*\tilde{\vh}|\}}\big)\\
		&\quad+\E\big(\Re(\tilde{\vh}^*\va\va^*\vx)\cdot I_{\{|\va^*\vx|\leq|\va^*\tilde{\vh}|\}}\big)\\
		&=2\cdot\E\big(\Re(\tilde{\vh}^*\va\va^*\vx)\cdot I_{\{|\va^*\vx|>|\va^*\tilde{\vh}|\}}\big),
		\end{align*}
		which implies (\ref{exp-haax}).
		
		Similarly, we have
		\begin{align*}
		&\E\big(|\va^*\vx|^2\big)\\
		&=\E\big(|\va^*\vx|^2\cdot I_{\{|\va^*\vx|>|\va^*\tilde{\vh}|\}}\big)+\E\big(|\va^*\vx|^2\cdot I_{\{|\va^*\vx|\leq|\va^*\tilde{\vh}|\}}\big)\\
    	&=\E\big(|\va^*\vx|^2\cdot I_{\{|\va^*\vx|>|\va^*\tilde{\vh}|\}}\big)+\E\big(|\va^*\tilde{\vh}|^2\cdot I_{\{|\va^*\vx|>\va^*\tilde{\vh}|\}}\big),
		\end{align*}
		which implies
		\[
		\E\big(|\va^*\vx|^2\cdot I_{\{|\va^*\vx|>|\va^*\tilde{\vh}|\}}\big)\geq \frac{1}{2}\E\big(|\va^*\vx|^2\big)= \frac{1}{2},
		\]	
		\[
		\E\big(|\va^*\vx|^2\cdot I_{\{|\va^*\vx|\leq|\va^*\tilde{\vh}|\}}\big)\leq \frac{1}{2}.
		\]	
		And
		\begin{align*}
		&2\cdot\E\bigg(\frac{\big(\Re(\tilde{\vh}^*\va\va^*\vx)\big)^2}{|\va^*\vx|^2}\cdot I_{\{|\va^*\vx|\leq|\va^*\tilde{\vh}|\}}\bigg)\\&\geq\E\bigg(\frac{\big(\Re(\tilde{\vh}^*\va\va^*\vx)\big)^2}{|\va^*\vx|^2}\bigg)=\E\bigg(\frac{\big(\Re(\tilde{\vh}^*\va\va^*\vx)\big)^2}{|\va^*\vx|^2}\cdot I_{\{|\va^*\vx|>|\va^*\tilde{\vh}|\}}\bigg)\\
		&\quad\quad+\E\bigg(\frac{\big(\Re(\tilde{\vh}^*\va\va^*\vx)\big)^2}{|\va^*\vx|^2}\cdot I_{\{|\va^*\vx|\leq|\va^*\tilde{\vh}|\}}\bigg)\\
		&\geq 2\cdot\E\bigg(\frac{\big(\Re(\tilde{\vh}^*\va\va^*\vx)\big)^2}{|\va^*\vx|^2}\cdot I_{\{|\va^*\vx|>|\va^*\tilde{\vh}|\}}\bigg)
		\end{align*}
		implies
		\[
		\E\bigg(\frac{\big(\Re(\tilde{\vh}^*\va\va^*\vx)\big)^2}{|\va^*\vx|^2}\cdot I_{\{|\va^*\vx|>|\va^*\tilde{\vh}|\}}\bigg)\leq\frac{1}{4}+\frac{1}{4}\Re^2(\tilde{\vh}^*\vx)
		\]
		and
		\[
		\E\bigg(\frac{\big(\Re(\tilde{\vh}^*\va\va^*\vx)\big)^2}{|\va^*\vx|^2}\cdot I_{\{|\va^*\vx|\leq|\va^*\tilde{\vh}|\}}\bigg)\geq\frac{1}{4}+\frac{1}{4}\Re^2(\tilde{\vh}^*\vx).
		\]
		Then to prove the inequalities (\ref{exp-ax1}), (\ref{exp-ax2}), (\ref{exp-haax/ax1}) and (\ref{exp-haax/ax2}), it's sufficient to prove
		\begin{equation}\label{ax_c1}
		\E\big( |\va^*\vx|^2\cdot I_{\{|\va^*\vx|>|\va^*\tilde{\vh}|\}}\big)\leq 3/4,
		\end{equation}
		\begin{equation}
		\begin{aligned}\label{haax/ax_c1}
		 \frac{1}{8}+\frac{7}{32}\Re^2(\tilde{\vh}^*\vx)\leq\E\bigg(\frac{\big(\Re(\tilde{\vh}^*\va\va^*\vx)\big)^2}{|\va^*\vx|^2}\cdot I_{\{|\va^*\vx|>|
			\va^*\tilde{\vh}\}}\bigg).
		\end{aligned}
		\end{equation}
		Next, we commit to prove (\ref{ax_c1}) and (\ref{haax/ax_c1}). Firstly, we take polar coordinates transformation:
		$$ \left\{
		\begin{aligned}
		\xi_1 & =  r_1\cos\theta_1 \\
		\xi_2 & =  r_1\sin\theta_1\\
		\xi_3 & =  (r_2\cos\theta_2-\sigma r_1\cos\theta_1)/\sqrt{1-\sigma^2}\\
		\xi_4 & =  (r_2\sin\theta_2-\sigma r_1\sin\theta_1)/\sqrt{1-\sigma^2}
		\end{aligned}
		\right.
		$$
		with $ r_1,r_2\in(0,\infty) $, $ \theta_1, \theta_2\in[0, 2\pi)$.
		Then we can write the expectation as
		\begin{align*}
		&\E\big( |\va^*\vx|^2\cdot I_{\{|\va^*\vx|>|\va^*\tilde{\vh}|\}}\big)\\
		 &=\frac{1}{\pi^2}\int_{0}^{2\pi}\int_{0}^{2\pi}\int_{0}^{\infty}\int_{0}^{r_1}\frac{r_1^3r_2}{1-\sigma^2}\exp\Big(-\frac{r_1^2+r_2^2}{1-\sigma^2}\Big)\\
		&\quad\cdot\exp\Big(\frac{2\sigma r_1r_2 \cos(\theta_1-\theta_2)}{1-\sigma^2}\Big)dr_2dr_1d\theta_1d\theta_2\\
		 &=4\sum_{k=0}^{\infty}\frac{1}{(k!)^2}\int_{0}^{\infty}\int_{0}^{r_1}\frac{\sigma^{2k}r_1^{2k+3}r_2^{2k+1}}{(1-\sigma^2)^{2k+1}}\\
		&\quad\quad \cdot\exp\Big(-\frac{r_1^2+r_2^2}{1-\sigma^2}\Big)dr_2dr_1\\
		&=\sum_{k=0}^{\infty}\sigma^{2k}(1-\sigma^2)^2(k+1)\Big(\frac{(2k+1)!!}{2^{k+2}\cdot k!}+\frac{k+1}{2}\Big).
		\end{align*}
		It is an even function about $ \sigma $ and when $ \sigma\in[0, 1) $ the derivative
		\begin{align*}
		&d\E\big( |\va^*\vx|^2\cdot I_{\{|\va^*\vx|>|\va^*\tilde{\vh}|\}}\big)/d\sigma\\
		%	     &\frac{d\E\big( |\va^*\vx|^2\cdot I_{\{|\va^*\vx|>|\va^*\tilde{\vh}|\}}\big)}{d\sigma}\\
		&=-\frac{\sigma}{4} - \sum_{k=1}^{\infty}\frac{(2k-1)!!}{2^{k+2}\cdot k!}\sigma^{2k+1}\leq 0.
		\end{align*}
		Hence the expectation obtains its maximum at $
		\sigma =0 $, i.e.,
		\begin{small}
			\begin{align*}
			&\E\big( |\va^*\vx|^2\cdot I_{\{|\va^*\vx|>|\va^*\tilde{\vh}|\}}\big)\\
			&\leq \frac{1}{\pi^2}\int_{0}^{2\pi}\int_{0}^{2\pi}\int_{0}^{\infty}\int_{0}^{r_1}r_1^3r_2\exp\big(-(r_1^2+r_2^2)\big)dr_2dr_1d\theta_1d\theta_2\\
			&=3/4.
			\end{align*}
		\end{small}
		Thus we have the inequality (\ref{ax_c1}).
		
		Using the same polar coordinates transformation,  we know
		\begin{align*}
		&\E\bigg(\frac{\big(\Re(\tilde{\vh}^*\va\va^*\vx)\big)^2}{|\va^*\vx|^2}\cdot I_{\{|\va^*\vx|>|
			\va^*\tilde{\vh}\}}\bigg)\\
		 &=\frac{1}{\pi^2}\int_{0}^{2\pi}\int_{0}^{2\pi}\int_{0}^{\infty}\int_{0}^{r_1}\frac{r_1r_2^3}{1-\sigma^2}\exp\Big(-\frac{r_1^2+r_2^2}{1-\sigma^2}\Big)\\
		&\quad\cdot\cos^2(\theta_1-\theta_2)\exp\Big(\frac{2\sigma r_1r_2 \cos(\theta_1-\theta_2)}{1-\sigma^2}\Big)dr_2dr_1d\theta_1d\theta_2\\
		 &=4\sum_{k=0}^{\infty}\frac{1}{(k!)^2}\frac{2k+1}{2k+2}\int_{0}^{\infty}\int_{0}^{r_1}\frac{\sigma^{2k}r_1^{2k+1}r_2^{2k+3}}{(1-\sigma^2)^{2k+1}}\\
		&\quad\quad \cdot\exp\Big(-\frac{r_1^2+r_2^2}{1-\sigma^2}\Big)dr_2dr_1\\
		 &=\frac{1}{2}\sum_{k=0}^{\infty}\sigma^{2k}(1-\sigma^2)^2\frac{1}{k!}\frac{2k+1}{2k+2}\Big((k+1)!-\frac{(2k+1)!!}{2^{k+1}} \Big)\\
		 &=\frac{1}{8}+\frac{7}{32}\sigma^2+\sum_{k=1}^{\infty}\frac{(2k+1)!!(2k+7)}{2^{k+4}(k+2)!}\sigma^{2k+2}.
		\end{align*}
		Thus we obtain (\ref{haax/ax_c1}). This completes the proof.
	\end{IEEEproof}
%%-----------------------------------------------------------------------------------------------------------
\begin{lemma}\label{concentration}
	Let $ \va_1, \va_2,\ldots,\va_m \in \C^n$  be i.i.d. Gaussian random measurements. Let $ \vx\in\C^n $ and $ \tilde{\vh}\in\C^n $ be two fixed vectors with $ \|\vx\|=\|\tilde{\vh}\|=1 $, $ \Im(\tilde{\vh}^*\vx)=0 $ and $ \tilde{\vh}\neq\pm \vx $. For any $ \delta>0 $, there exist positive constants $C_\delta, c_\delta>0$ such that for any $m \geq C_\delta n$ the inequalities
	\begin{align}\label{con1}
	\bigg|\frac{1}{m}\sum_{j=1}^{m}	 \Re(\tilde{\vh}^*\va_j\va_j^*\vx)\cdot I_{\{|\va_j^*\vx|>|\va_j^*\tilde{\vh}|\}}-\frac{1}{2}\Re(\tilde{\vh}^*\vx)\bigg|	 \leq \delta,
	\end{align}
	\begin{align}\label{con2}
	\frac{1}{2}-\delta\leq \frac{1}{m}\sum_{j=1}^{m}|\va_j^*\vx|^2\cdot I_{\{|\va_j^*\vx|>|\va_j^*\tilde{\vh}|\}}\leq \frac{3}{4}+\delta,
	\end{align}
		\begin{align}\label{con2'}
	\frac{1}{4}-\delta\leq \frac{1}{m}\sum_{j=1}^{m}|\va_j^*\vx|^2\cdot I_{\{|\va_j^*\vx|\leq|\va_j^*\tilde{\vh}|\}}\leq \frac{1}{2}+\delta,
	\end{align}
	\begin{equation}
	\begin{aligned}\label{con3}
	&\frac{1}{8}+\frac{7}{32}\Re^2(\tilde{\vh}^*\vx)-\delta\\
	&\quad\leq\frac{1}{m}\sum_{j=1}^{m}	 \frac{\big(\Re(\tilde{\vh}^*\va_j\va_j^*\vx)\big)^2}{|\va_j^*\vx|^2}\cdot I_{\{|\va_j^*\vx|>|\va_j^*\tilde{\vh}|\}}\\
	&\quad\leq\frac{1}{4}+\frac{1}{4}\Re^2(\tilde{\vh}^*\vx)+\delta
	\end{aligned}
	\end{equation}
	and
	\begin{equation}
	\begin{aligned}\label{con4}
	&\frac{1}{4}+\frac{1}{4}\Re^2(\tilde{\vh}^*\vx)-\delta\\
	&\quad\leq\frac{1}{m}\sum_{j=1}^{m}	 \frac{\big(\Re(\tilde{\vh}^*\va_j\va_j^*\vx)\big)^2}{|\va_j^*\vx|^2})\cdot I_{\{|\va_j^*\vx|\leq|\va_j^*\tilde{\vh}|\}}\\
	&\quad\leq\frac{3}{8}+\frac{9}{32}\Re^2(\tilde{\vh}^*\vx)+\delta
	\end{aligned}
	\end{equation}
	hold with probability at least $ 1-\exp(-c_\delta m) $.
\end{lemma}
\begin{IEEEproof}
	For fixed $ \tilde{\vh} $ and $ \vx $, the following sets are all independent sub-exponential random variables
	$$  \big\{ \Re(\tilde{\vh}^*\va_j\va_j^*\vx)\cdot I_{\{|\va_j^*\vx|>|\va_j^*\tilde{\vh}|\}},\, j=1,\ldots,m\big\}, $$
	$$  \big\{ |\va_j^*\vx|^2\cdot I_{\{|\va_j^*\vx|>|\va_j^*\tilde{\vh}|\}},\, j=1,\ldots,m\big\}, $$
	$$ \big \{ |\va_j^*\vx|^2\cdot I_{\{|\va_j^*\vx|\leq|\va_j^*\tilde{\vh}|\}},\, j=1,\ldots,m\big\}, $$
	$$  \Big\{ \frac{\big(\Re(\tilde{\vh}^*\va_j\va_j^*\vx)\big)^2}{|\va_j^*\vx|^2}\cdot I_{\{|\va_j^*\vx|>|\va_j^*\tilde{\vh}|\}},\, j=1,\ldots,m\Big\}, $$
	$$  \Big\{ \frac{\big(\Re(\tilde{\vh}^*\va_j\va_j^*\vx)\big)^2}{|\va_j^*\vx|^2}\cdot I_{\{|\va_j^*\vx|\leq|\va_j^*\tilde{\vh}|\}},\, j=1,\ldots,m\Big\}. $$
	Recall that $ \va=(a_1,\ldots,a_n)\in\C^n \sim \mathcal{N}(0,I/2) + i\mathcal{N}(0,I/2) $ is a Gaussian random measurement.
	Then based on Bernstein-type inequality, for any $ \delta >0$, the inequalities
	\begin{align*}
	\bigg|\frac{1}{m}\sum_{j=1}^{m}&	 \Re(\tilde{\vh}^*\va_j\va_j^*\vx)\cdot I_{\{|\va_j^*\vx|>|\va_j^*\tilde{\vh}|\}}\\
	&-\E\big(\Re(\tilde{\vh}^*\va\va^*\vx)\cdot I_{\{|\va*\vx|>|\va^*\tilde{\vh}|\}}\big)\bigg|	 \leq \delta,
	\end{align*}
	\begin{align*}
	\bigg|\frac{1}{m}\sum_{j=1}^{m}|\va_j^*\vx|^2\cdot I_{\{|\va_j^*\vx|>|\va_j^*\tilde{\vh}|\}}-\E(|\va^*\vx|^2\cdot I_{\{|\va^*\vx|>|
		\va^*\tilde{\vh}\}})\bigg|\leq\delta,
	\end{align*}
	\begin{align*}
	\bigg|\frac{1}{m}\sum_{j=1}^{m}|\va_j^*\vx|^2\cdot I_{\{|\va_j^*\vx|\leq|\va_j^*\tilde{\vh}|\}}-\E(|\va^*\vx|^2\cdot I_{\{|\va^*\vx|\leq|
		\va^*\tilde{\vh}\}})\bigg|\leq\delta,
	\end{align*}
	\begin{align*}
	\bigg|\frac{1}{m}\sum_{j=1}^{m}	 &\frac{\big(\Re(\tilde{\vh}^*\va_j\va_j^*\vx)\big)^2}{|\va_j^*\vx|^2}\cdot I_{\{|\va_j^*\vx|>|\va_j^*\tilde{\vh}|\}}\\
	&-\E\bigg(\frac{\big(\Re(\tilde{\vh}^*\va\va^*\vx)\big)^2}{|\va^*\vx|^2}\cdot I_{\{|\va_j^*\vx|>|\va_j^*\tilde{\vh}|\}}\bigg)\bigg|\leq\delta,
	\end{align*}
	\begin{align*}
	\bigg|\frac{1}{m}\sum_{j=1}^{m}	 &\frac{\big(\Re(\tilde{\vh}^*\va_j\va_j^*\vx)\big)^2}{|\va_j^*\vx|^2}\cdot I_{\{|\va_j^*\vx|\leq|\va_j^*\tilde{\vh}|\}}\\
	&-\E\bigg(\frac{\big(\Re(\tilde{\vh}^*\va\va^*\vx)\big)^2}{|\va^*\vx|^2}\cdot I_{\{|\va_j^*\vx|\leq|\va_j^*\tilde{\vh}|\}}\bigg)\bigg|\leq\delta
	\end{align*}
	hold with probability at least $ 1-\exp(-c_\delta m)  $ provided $ m\geq C_\delta n $, where $C_\delta, c_\delta$ are positive constants depending on $\delta$.  Then the inequalities (\ref{con1}), (\ref{con2}), (\ref{con2'}), (\ref{con3}), (\ref{con4}) can be derived directly from the expectation bounds given in Lemma \ref{expectation}.

\end{IEEEproof}
%%-----------------------------------------------------------------------------------------------------------

	The following lemma provides an upper bound for the operator norm of $  \nabla^2 f_\vepsilon(\vz)$.

	\begin{lemma}\label{hessian}
		Set $  \vepsilon = \sqrt{\alpha}\vb  $. Then there exist constants $C',c'>0$ such that for $ m\geq C'n $,
		$\|\nabla^2 f_\vepsilon(\vz)\|\leq2\sqrt{\frac{1+\alpha}{\alpha}}  $ holds with probability at least $ 1-\exp(-c'm) $.
	\end{lemma}
	\begin{IEEEproof}
		Recall that
		\begin{align*}
		\nabla f_{\vepsilon}(\vz) &:= \left(\frac{\partial f_{\vepsilon}(\vz,\overline{\vz})}{\partial \vz}\Big|_{\overline{\vz} = \text{constant}}\right)^* \\
		& = \frac{1}{m}\sum_{j=1}^{m}\left(1-\frac{\sqrt{b_j^2 +\epsilon_j^2}}{\sqrt{|\va_j^*\vz|^2 +\epsilon_j^2}}\right)\va_j\va_j^*\vz.
		\end{align*}
		Similarly,  we obtain
		\begin{align*}
		&\nabla^2 f_{\vepsilon}(\vz)\\
		&=  \frac{1}{m}\sum_{j=1}^{m}\left(1-\frac{\sqrt{b_j^2 +\epsilon_j^2}}{\sqrt{|\va_j^*\vz|^2 +\epsilon_j^2}}\right)\va_j\va_j^*\\
		& \quad\quad +\frac{1}{m}\sum_{j=1}^{m}
		 \Bigg(\frac{\sqrt{b_j^2+\epsilon_j^2}|\va_j^*\vz|^2}{2\big(|\va_j^*\vz|^2+\epsilon_j^2\big)^{3/2}}\Bigg)\va_j\va_j^*\\
		 &=\frac{1}{m}\sum_{j=1}^{m}\left(1-\frac{\sqrt{b_j^2+\epsilon_j^2}\big(\frac{1}{2}|\va_j^*\vz|^2+\epsilon_j^2\big)}{\big(|\va_j^*\vz|^2+\epsilon_j^2\big)^{3/2}}\right)\va_j\va_j^*.\\
		\end{align*}
		For any $ \vz \in \C^n$, we have
		\begin{align*}
		&\|\nabla^2 f_{\vepsilon}(\vz)\|\\
		 &=\max_{\vy\in\SS^{n-1}}\frac{1}{m}\sum_{j=1}^{m}\left(1-\frac{\sqrt{b_j^2+\epsilon_j^2}\big(\frac{1}{2}|\va_j^*\vz|^2+\epsilon_j^2\big)}{\big(|\va_j^*\vz|^2+\epsilon_j^2\big)^{3/2}}\right)|\va_j^*\vy|^2\\
		&\leq \max_{\vy\in\SS^{n-1}}\frac{1}{m}\sum_{j=1}^{m}\left(1+\frac{\sqrt{b_j^2+\epsilon_j^2}}{\sqrt{|\va_j^*\vz|^2+\epsilon_j^2}}\right)|\va_j^*\vy|^2\\
		&\leq \max_{\vy\in\SS^{n-1}}\frac{1}{m}\sum_{j=1}^{m}\left(1+\sqrt{\frac{1+\alpha}{\alpha}}\right)
		|\va_j^*\vy|^2\\
		&\leq 2\sqrt{\frac{1+\alpha}{\alpha}}
		\end{align*}
		with probability at least $ 1-\exp(-c'm) $ provided $ m\geq C'n $. Here the third inequality is obtained by Lemma \ref{sub_gaussian_concentration}.
	\end{IEEEproof}

	% =========================Sec===============================

	% =========================Ref===============================
	%\bibliographystyle{IEEEtran}
	\bibliographystyle{IEEEtran}
	\bibliography{phase_retrieval_newmodel}

\begin{thebibliography}{10}
\providecommand{\url}[1]{#1}
\csname url@samestyle\endcsname
\providecommand{\newblock}{\relax}
\providecommand{\bibinfo}[2]{#2}
\providecommand{\BIBentrySTDinterwordspacing}{\spaceskip=0pt\relax}
\providecommand{\BIBentryALTinterwordstretchfactor}{4}
\providecommand{\BIBentryALTinterwordspacing}{\spaceskip=\fontdimen2\font plus
\BIBentryALTinterwordstretchfactor\fontdimen3\font minus
  \fontdimen4\font\relax}
\providecommand{\BIBforeignlanguage}[2]{{%
\expandafter\ifx\csname l@#1\endcsname\relax
\typeout{** WARNING: IEEEtran.bst: No hyphenation pattern has been}%
\typeout{** loaded for the language `#1'. Using the pattern for}%
\typeout{** the default language instead.}%
\else
\language=\csname l@#1\endcsname
\fi
#2}}
\providecommand{\BIBdecl}{\relax}
\BIBdecl

\bibitem{miao1999extending}
J.~Miao, P.~Charalambous, J.~Kirz, and D.~Sayre, ``Extending the methodology of
  x-ray crystallography to allow imaging of micrometre-sized non-crystalline
  specimens,'' \emph{Nature}, vol. 400, no. 6742, p. 342, 1999.

\bibitem{elser2018benchmark}
V.~Elser, T.~Lan, and T.~Bendory, ``Benchmark problems for phase retrieval,''
  \emph{Siam Journal on Imaging Sciences}, vol.~11, no.~4, pp. 2429--2455,
  2018.

\bibitem{fienup1987phase}
C.~Fienup and J.~Dainty, ``Phase retrieval and image reconstruction for
  astronomy,'' \emph{Image Recovery: Theory and Application}, pp. 231--275,
  1987.

\bibitem{walther1963question}
A.~Walther, ``The question of phase retrieval in optics,'' \emph{Optica Acta:
  International Journal of Optics}, vol.~10, no.~1, pp. 41--49, 1963.

\bibitem{millane1990phase}
R.~P. Millane, ``Phase retrieval in crystallography and optics,'' \emph{JOSA
  A}, vol.~7, no.~3, pp. 394--411, 1990.

\bibitem{miao2008extending}
J.~Miao, T.~Ishikawa, Q.~Shen, and T.~Earnest, ``Extending x-ray
  crystallography to allow the imaging of noncrystalline materials, cells, and
  single protein complexes,'' \emph{Annu. Rev. Phys. Chem.}, vol.~59, pp.
  387--410, 2008.

\bibitem{balan2006signal}
R.~Balan, P.~Casazza, and D.~Edidin, ``On signal reconstruction without
  phase,'' \emph{Applied and Computational Harmonic Analysis}, vol.~20, no.~3,
  pp. 345--356, 2006.

\bibitem{bandeira2014saving}
A.~S. Bandeira, J.~Cahill, D.~G. Mixon, and A.~A. Nelson, ``Saving phase:
  Injectivity and stability for phase retrieval,'' \emph{Applied and
  Computational Harmonic Analysis}, vol.~37, no.~1, pp. 106--125, 2014.

\bibitem{WX17}
Y.~Wang and Z.~Xu, ``Generalized phase retrieval: Measurement number, matrix
  recovery and beyond,'' \emph{Applied and Computational Harmonic Analysis},
  vol.~47, no.~2, pp. 423--446, 2017.

\bibitem{huang2016phase}
K.~Huang, Y.~C. Eldar, and N.~D. Sidiropoulos, ``Phase retrieval from 1d
  fourier measurements: Convexity, uniqueness, and algorithms,'' \emph{IEEE
  Transactions on Signal Processing}, vol.~64, no.~23, pp. 6105--6117, 2016.

\bibitem{eldar2015sparse}
Y.~C. Eldar, P.~Sidorenko, D.~G. Mixon, S.~Barel, and O.~Cohen, ``Sparse phase
  retrieval from short-time fourier measurements,'' \emph{IEEE Signal
  Processing Letters}, vol.~22, no.~5, pp. 638--642, 2015.

\bibitem{jaganathan2016stft}
K.~Jaganathan, Y.~C. Eldar, and B.~Hassibi, ``Stft phase retrieval: Uniqueness
  guarantees and recovery algorithms,'' \emph{IEEE Journal of Selected Topics
  in Signal Processing}, vol.~10, no.~4, pp. 770--781, 2016.

\bibitem{bendory2018non-convex}
T.~Bendory, Y.~C. Eldar, and N.~Boumal, ``Non-convex phase retrieval from stft
  measurements,'' \emph{IEEE Transactions on Information Theory}, vol.~64,
  no.~1, pp. 467--484, 2018.

\bibitem{candes2013phaselift}
E.~J. Cand{\`e}s, T.~Strohmer, and V.~Voroninski, ``Phaselift: Exact and stable
  signal recovery from magnitude measurements via convex programming,''
  \emph{Communications on Pure and Applied Mathematics}, vol.~66, no.~8, pp.
  1241--1274, 2013.

\bibitem{candes_phaselift}
E.~J. Cand{\`e}s and X.~Li, ``Solving quadratic equations via phaselift when
  there are about as many equations as unknowns,'' \emph{Foundations of
  Computational Mathematics}, vol.~14, no.~5, pp. 1017--1026, 2014.

\bibitem{candes_matrix_completion}
E.~J. Cand{\`e}s, Y.~C. Eldar, T.~Strohmer, and V.~Voroninski, ``Phase
  retrieval via matrix completion,'' \emph{SIAM review}, vol.~57, no.~2, pp.
  225--251, 2015.

\bibitem{goldstein2018phasemax}
T.~Goldstein and C.~Studer, ``Phasemax: Convex phase retrieval via basis
  pursuit,'' \emph{IEEE Transactions on Information Theory}, 2018.

\bibitem{dhifallah2017fundamental}
O.~Dhifallah and Y.~M. Lu, ``Fundamental limits of phasemax for phase
  retrieval: A replica analysis,'' in \emph{Computational Advances in
  Multi-Sensor Adaptive Processing (CAMSAP), 2017 IEEE 7th International
  Workshop on}.\hskip 1em plus 0.5em minus 0.4em\relax IEEE, 2017, pp. 1--5.

\bibitem{hand2016elementary}
P.~Hand and V.~Voroninski, ``An elementary proof of convex phase retrieval in
  the natural parameter space via the linear program phasemax,''
  \emph{Communications in Mathematical Sciences}, vol.~16, no.~7, pp.
  2047--2051, 2018.

\bibitem{gerchberg1972practical}
R.~W. Gerchberg, ``A practical algorithm for the determination of phase from
  image and diffraction plane pictures,'' \emph{Optik}, vol.~35, pp. 237--246,
  1972.

\bibitem{fienup1982phase}
J.~R. Fienup, ``Phase retrieval algorithms: a comparison,'' \emph{Applied
  optics}, vol.~21, no.~15, pp. 2758--2769, 1982.

\bibitem{netrapalli2013phase}
P.~Netrapalli, P.~Jain, and S.~Sanghavi, ``Phase retrieval using alternating
  minimization,'' in \emph{Advances in Neural Information Processing Systems},
  2013, pp. 2796--2804.

\bibitem{waldspurger2018phase}
I.~Waldspurger, ``Phase retrieval with random gaussian sensing vectors by
  alternating projections,'' \emph{IEEE Transactions on Information Theory},
  vol.~64, no.~5, pp. 3301--3312, 2018.

\bibitem{candes_wf}
E.~J. Cand{\`e}s, X.~Li, and M.~Soltanolkotabi, ``Phase retrieval via wirtinger
  flow: Theory and algorithms,'' \emph{Information Theory, IEEE Transactions
  on}, vol.~61, no.~4, pp. 1985--2007, 2015.

\bibitem{gao2017phaseless}
B.~Gao and Z.~Xu, ``Phaseless recovery using the gauss-newton method,''
  \emph{IEEE Transactions on Signal Processing}, vol.~65, no.~22, pp.
  5885--5896, 2017.

\bibitem{wang2018solving}
G.~Wang, G.~B. Giannakis, and Y.~C. Eldar, ``Solving systems of random
  quadratic equations via truncated amplitude flow,'' \emph{IEEE Transactions
  on Information Theory}, vol.~64, no.~2, pp. 773--794, 2018.

\bibitem{zhang2016reshaped}
H.~Zhang and Y.~Liang, ``Reshaped wirtinger flow for solving quadratic system
  of equations,'' in \emph{Advances in Neural Information Processing Systems},
  2016, pp. 2622--2630.

\bibitem{huangxu}
M.~Huang and Z.~Xu, ``The estimation performance of nonlinear least squares for
  phase retrieval,'' \emph{IEEE Transactions on Information Theory}, 2020.

\bibitem{wang2018phase}
G.~Wang, G.~B. Giannakis, Y.~Saad, and J.~Chen, ``Phase retrieval via
  reweighted amplitude flow,'' \emph{IEEE Transactions on Signal Processing},
  vol.~66, no.~11, pp. 2818--2833, 2018.

\bibitem{chen2015solving}
Y.~Chen and E.~Cand{\`e}s, ``Solving random quadratic systems of equations is
  nearly as easy as solving linear systems,'' in \emph{Communications on Pure
  and Applied Mathematics}, vol.~70, 2017, pp. 822--883.

\bibitem{sun2018geometric}
J.~Sun, Q.~Qu, and J.~Wright, ``A geometric analysis of phase retrieval,''
  \emph{Foundations of Computational Mathematics}, vol.~18, no.~5, pp.
  1131--1198, 2018.

\bibitem{2017Implicit}
C.~Ma, K.~Wang, Y.~Chi, and Y.~Chen, ``Implicit regularization in nonconvex
  statistical estimation: Gradient descent converges linearly for phase
  retrieval, matrix completion and blind deconvolution,'' \emph{Foundations of
  Computational Mathematics}, vol.~20, no.~3, pp. 451--632, 2020.

\bibitem{Zhang2017A}
H.~Zhang, Y.~Zhou, Y.~Liang, and Y.~Chi, ``A nonconvex approach for phase
  retrieval: reshaped wirtinger flow and incremental algorithms,''
  \emph{Journal of Machine Learning Research}, vol.~18, no. 141, pp. 1--35,
  2017.

\bibitem{cdp}
E.~J. Candes, X.~Li, and M.~Soltanolkotabi, ``Phase retrieval from coded
  diffraction patterns,'' \emph{Applied and Computational Harmonic Analysis},
  vol.~39, no.~2, pp. 277--299, 2015.

\bibitem{vershynin2010introduction}
R.~Vershynin, ``Introduction to the non-asymptotic analysis of random
  matrices,'' \emph{Compressed Sensing, Theory and Applications}, 2012.

\end{thebibliography}
\end{document}